\documentclass[hidelinks,onefignum,onetabnum]{siamart220329}

\usepackage{lipsum}
\usepackage{amsfonts,amssymb,booktabs}
\usepackage{graphicx}
\usepackage{epstopdf}
\usepackage{algorithmic}
\ifpdf
  \DeclareGraphicsExtensions{.eps,.pdf,.png,.jpg}
\else
  \DeclareGraphicsExtensions{.eps}
\fi

\newsiamremark{remark}{Remark}
\newsiamremark{hypothesis}{Hypothesis}
\crefname{hypothesis}{Hypothesis}{Hypotheses}
\newsiamthm{claim}{Claim}

\headers{Analytical Galerkin layer potential quadrature}{ N. A. Gumerov, S. Kaneko, and R. Duraiswami}

\title{Analytical Galerkin boundary integrals of Laplace kernel layer potentials in $\mathbb{R}^3$ \thanks{Submitted to the editors December 6th, 2022.
\funding{Cooperative Research Agreements W911NF1420118 and W911NF2020213 between UMD and Army Research Laboratory, with D. Hull, R. Adelman and S. Vinci as Technical monitors. S.~Kaneko thanks Japan Student Services Organization and Watanabe Foundation for scholarships.}}}

\author{
Nail A. Gumerov\thanks{UMIACS, University of Maryland, College Park. Author passed away during paper preparation.}
\and 
Shoken Kaneko\thanks{Department of Computer Science and UMIACS, University of Maryland, College Park. (\email{kaneko60@umd.edu}, \email{ramanid@umd.edu}).  } 
\and
Ramani Duraiswami \footnotemark[3]
}

\usepackage{amsopn}

\ifpdf
\hypersetup{
  pdftitle={An Example Article},
  pdfauthor={D. Doe, P. T. Frank, and J. E. Smith}
}
\fi

\allowdisplaybreaks

\begin{document}
\raggedbottom
\maketitle
\begin{abstract}
A method for analytical computation of the double surface integrals for all layer potential kernels 
associated with the Laplace Green's function, in the Galerkin boundary element method (BEM) 
in $\mathbb{R}^3$ using piecewise constant flat elements is presented. 
The method uses recursive dimensionality reduction from 4D  ($\mathbb{R}^2\times\mathbb{R}^2$) based on Gauss' divergence theorem. 
Computable analytical expressions for all cases of relative location of the source and receiver triangles are covered 
for the single and double layer potentials and their gradients with analytical treatment of the singular cases are presented. 
A trick that enables reduction of the case of gradient of the single layer to the same integrals as for the single layer is introduced 
using symmetry properties. 
The method was confirmed using analytical benchmark cases, comparisons with error-controlled computations of regular multidimensional integrals, and a convergence study for singular cases.
\end{abstract}



\section{Introduction}

The Galerkin BEM is a commonly used variant of general method of moments (MoM), 
where the basis and test functions are the same \cite{Harrington1983book}, and is often preferred as it yields symmetric matrices. 
A challenge with the method is efficient and  accurate evaluation of double surface ($\mathbb{R}^2\times\mathbb{R}^2$) integrals.  
When the source and receiver triangles are well separated (regular integrals) efficient numerical procedures using exact 
analytical expressions for the potentials over the source triangles \cite{Hess1967PAS, Salvadori2001CMAME, Gumerov2021arXiv} and 
Gaussian  quadrature over the receiver triangle can be employed. 
Moreover, the Galerkin integrals can be computed using multipole expansions \cite{Adelman2016IEEE, Gumerov2021arXiv2}. 
When the triangles have one or more common points all integrals are singular, and treated using special methods 
(e.g., \cite{Adelman2016IEEE, DEllia2011IJNMBE}). 
With two common vertices between source and receiver triangles the integrals are hypersingular and do not exist even in the 
Cauchy principal value sense, and requires a special treatment \cite{Salvadori2001CMAME,Gray2001EABE,Polimeridis2011IEEE}. 
Depending on triangle proximity we also have the case of nearly singular integrals.

While the literature provides numerical methods for the Galerkin integrals,  full analytical expressions  were derived and analyzed only recently.  
Previously analytical results were only available for the self-action Galerkin integrals \cite{Eibert1995IEEE, Sievers2005IEEE}. 
The ``modern'' history of analytical expressions for the Galerkin integrals probably starts with \cite{Rosen1993SIAM}, 
who  used Gauss' divergence theorem to reduce the volume integrals to surface integrals for homogeneous functions and introduced the concept of the primitive boundary function (PBF), 
though not in the context of the Galerkin method. The PBF is the integrand of the boundary integral that is related to the original integrand of the volume integral. 

Lenoir \cite{Lenoir2006MASP} realized then that the PBF can be computed recursively for  $d=4,3,2,1$ to reduce the 4D integrals for the single layer potential to a final form, 
which contains only elementary functions. This was restricted to the case when the planes of the source and receiver triangles are not parallel. 
Next Lenoir and Salles \cite{Lenoir2012SIAM} considered the special case when source and receiver triangles are located in the same plane or parallel planes. 
The gradient of single layer is calculated for this special case. 
The same authors combined these two results and added the case of linear basis functions in \cite{Lenoir2012ECCOMAS, Salles2013Thesis}. 
In \cite{Warncke2019SMAI} the expressions were further simplified and the number of cases needed to be considered reduced. 
An expression for the hypersingular integral was also presented, though no treatment of the singularity was provided. 
A related recent study is \cite{Oueslati2021IJNME}, where the focus is on  single and double layer integrals. 
Analytical expressions for the double layer were obtained only when the source and receiver triangles have a common vertex or edge.

In comparison to previous work, we follow an independent approach to recursive dimensionality reduction, 
and a different triangle parametrization and treatment of the integrand. In our work this appears generally 
as a 4-parametric non-homogeneous function with 5 three dimensional vector coefficients. 
Homogenization of the PBF occurs via projections which can be efficiently generally handled using Gram-Schmidt orthogonalization. 
Further, we provide explicit expressions for the PBF in all dimensions and we found substantial simplifications related to the 
interplay of the 3D real space and 4D space in which integration occurs. 
We consider all cases of relative locations of the source and receiver triangles, and all layer potentials - including single and 
double layer potentials, the gradient of single layer potential, and the hypersingular case with analytical treatment of the singularity. 
This provides analytical expressions for the cases which were not considered before. 
We were able to reduce the gradient of the single layer to the same integrals as for the single layer using symmetry properties, 
allowing amortization of the cost of computation of the gradient and also the double layer potential.  
We confirmed our expressions via analytical benchmark cases, comparisons with error-controlled computations of 
regular multidimensional integrals and a convergence study for singular cases. 
For reasons of space, example illustrations for self-action integral case are provided in the supplemental material.

\section{Problem Statement}
We seek analytical expressions for:
\begin{eqnarray}
L &=&\int_{S_{y}}\int_{S_{x}}G\left( \mathbf{x,y}\right) dS\left( \mathbf{x}%
\right) dS\left( \mathbf{y}\right) ,  \label{ps1} \\
\text{ }M &=&\int_{S_{y}}\int_{S_{x}}\mathbf{n}_{x}\cdot \nabla _{\mathbf{x}%
}G\left( \mathbf{x,y}\right) dS\left( \mathbf{x}\right) dS\left( \mathbf{y}%
\right) ,  \label{ps2} \\
\mathbf{L}^{\prime } &=&\int_{S_{y}}\nabla _{\mathbf{y}}\int_{S_{x}}G\left( 
\mathbf{x,y}\right) dS\left( \mathbf{x}\right) dS\left( \mathbf{y}\right) ,
\label{ps3} \\
M^{\prime } &=&\int_{S_{y}}\mathbf{n}_{y}\cdot \nabla _{\mathbf{y}%
}\int_{S_{x}}\mathbf{n}_{x}\cdot \nabla _{\mathbf{x}}G\left( \mathbf{x,y}%
\right) dS\left( \mathbf{x}\right) dS\left( \mathbf{y}\right) ,  \label{ps4}
\end{eqnarray}%
where $S_{x}$ and $S_{y}$ are arbitrary positively oriented flat triangles
in $\mathbb{R}^3$ with normals $\mathbf{n}_{x}$ and $\mathbf{n}_{y}$, and vertices 
$\left( \mathbf{x}_{1},\mathbf{x}_{2},\mathbf{x}_{3}\right) $ and $\left( 
\mathbf{y}_{1},\mathbf{y}_{2},\mathbf{y}_{3}\right) $; while\ 
$\nabla _{x}$ and $\nabla _{y}$ are the gradient operators with respect to $\mathbf{x}$ and $\mathbf{y}$; 
and $G\left( \mathbf{x,y}\right) $ the free-space Green's function for the Laplace equation is given by
\begin{equation}
G\left( \mathbf{x,y}\right) =\frac{1}{\left| \mathbf{x}-\mathbf{y}\right| }.
\label{ps6}
\end{equation}

Integrals (\ref{ps1}) and (\ref{ps2}) are the single and double layer potentials, (\ref{ps3}) the gradient of single layer 
used sometimes for reduction of hypersingularity or in boundary integral formulations involving surface gradients or curls 
(electromagnetism), and (\ref{ps4}) the integral of the normal derivative of the double layer potential. 
All integrals are singular. In case the singularity is not integrable, the integral is understood in terms of the principal value. 
There are special cases when integral $M^{\prime }$ diverges even in the sense of principal value 
(when $S_{x}$ and $S_{y}$ have a common edge), for which we provide a special consideration. 

\section{Single layer potential}

\subsection{Standard triangle}

Using a simple linear transformation, integrals over arbitrary triangles $S_{x}$ and $S_{y}$ can be transformed to 
integration over standard triangles $T_{x}$ and $T_{y}$, with vertices $\left( 0,0\right) ,\left( 1,0\right) ,$ and $\left(0,1\right) $. 
Integrals over the source and target triangles respectively have the (independent) parametrizations,%
\begin{gather}
\mathbf{x=X}\left( s,t\right) ,\quad \mathbf{y=Y}\left( u,v\right) .
\label{s1}\\
\mathbf{X}\left( s,t\right) = \mathbf{X}_{s}s+\mathbf{X}_{t}t+\mathbf{X}%
_{0},  \label{s2} \quad
\mathbf{Y}\left( u,v\right) = \mathbf{Y}_{u}u+\mathbf{Y}_{v}v+\mathbf{Y}%
_{0}. 
\end{gather}%
The vector constants are uniquely determined by the vertices. 
Using the  correspondence of the vertices with indices 1, 2, and 3 to $\left( 0,0\right) ,\left( 1,0\right) ,$ and $\left( 0,1\right) $, we obtain
\begin{equation}\begin{aligned}
\mathbf{x}_{1} &= \mathbf{X}_{0},\quad \mathbf{x}_{2}=\mathbf{X}_{s}+\mathbf{%
X}_{0},\quad \mathbf{x}_{3}=\mathbf{X}_{t}+\mathbf{X}_{0},\\  \label{s3} \quad
\mathbf{y}_{1} &= \mathbf{Y}_{0},\quad \mathbf{y}_{2}=\mathbf{Y}_{u}+\mathbf{%
Y}_{0},\quad \mathbf{y}_{3}=\mathbf{Y}_{v}+\mathbf{Y}_{0}. 
\end{aligned}\end{equation}%
This can be resolved as%
\begin{equation}\begin{aligned}
\mathbf{X}_{0} &= \mathbf{x}_{1},\quad \mathbf{X}_{s}=\mathbf{x}_{2}-\mathbf{%
x}_{1},\quad \mathbf{X}_{t}=\mathbf{x}_{3}-\mathbf{x}_{1},  \\\label{s4} \quad
\mathbf{Y}_{0} &= \mathbf{y}_{1},\quad \mathbf{Y}_{u}=\mathbf{y}_{2}-\mathbf{%
y}_{1},\quad \mathbf{Y}_{v}=\mathbf{y}_{3}-\mathbf{y}_{1}. 
\end{aligned}\end{equation}%
The Jacobians of the two transforms to standard triangles are respectively 
\begin{equation}
J_{x}=\left| \mathbf{X}_{s}\times \mathbf{X}_{t}\right| =2A_{x},\quad
J_{y}=\left| \mathbf{Y}_{u}\times \mathbf{Y}_{v}\right| =2A_{y},  \label{s5}
\end{equation}%
where $A_{x}$ and $A_{y}$ are the areas of $S_{x}$ and $S_{y}$. 
This transform shows that%
\begin{equation}\begin{aligned}
& L=\int_{S_{y}}\int_{S_{x}} G\left( \mathbf{x,y}\right) dS\left( \mathbf{x}%
\right) dS\left( \mathbf{y}\right) 
= 4A_{x}A_{y}U,  \label{s6} \quad 
g\left( s,t,u,v\right) = G\left( \mathbf{X}\left( s,t\right) \mathbf{,Y}%
\left( u,v\right) \right),\\
& U = \int_{T_{y}}\int_{T_{x}}g\left( s,t,u,v\right) dS\left( s,t\right)
dS\left( u,v\right)  = \int_{0}^{1}\int_{0}^{1-u}\int_{0}^{1}\int_{0}^{1-s}g\left(
s,t,u,v\right) dtdsdvdu. 
\end{aligned}\end{equation}
As needed, we will also use the following notation for triangle $S_{x}$%
\begin{equation}
\mathbf{l}_{x1}=\mathbf{x}_{2}-\mathbf{x}_{1},\quad \mathbf{l}_{x2}=\mathbf{x%
}_{3}-\mathbf{x}_{2},\quad \mathbf{l}_{x3}=\mathbf{x}_{1}-\mathbf{x}%
_{3},\quad l_{xj}=\left| \mathbf{l}_{xj}\right| ,\quad j=1,2,3,  \label{s10}
\end{equation}%
and with subscript $y$ for triangle $S_{y}$. Comparing  with (\ref{s4}), we find%
\begin{equation}\begin{aligned}
\mathbf{X}_{s} &= \mathbf{l}_{x1},\quad \mathbf{X}_{t}=-\mathbf{l}%
_{x3},\quad \mathbf{l}_{x2}=\mathbf{X}_{t}-\mathbf{X}_{s},  \label{s11} \\
\mathbf{Y}_{u} &= \mathbf{l}_{y1},\quad \mathbf{Y}_{v}=-\mathbf{l}%
_{y3},\quad \mathbf{l}_{y2}=\mathbf{Y}_{v}-\mathbf{Y}_{u}. 
\end{aligned}\end{equation}

\subsection{Multidimensional parametrization and integrals}
We introduce a recursive technique, which 
reduces integrals over $d$-dimensional manifolds to sums of integrals over $\left( d-1\right) $-dimensional manifolds. 
The method applies to any function $G\left( \left| \mathbf{x-y}\right|\right)$ dependent on the distance 
between the source and evaluation points, such as a Green's function. 
Let points in $\mathbb{R}^{d}$ be characterized by coordinates $\left( s_{1},...,s_{d}\right) $ corresponding to $d$-dimensional basis vectors $\mathbf{i}_{1},...,\mathbf{i}_{d}$. 
Let now $\mathbf{a}_{1}\mathbf{,...,a}_{d}\mathbf{,e}_{d}\in \mathbb{R}^{3}$ be some constant 3D vectors (parameters), which define a linear function
\begin{equation}
\mathbf{R}_{d}\left( s_{1},...,s_{d}\right) =\mathbf{a}_{1}s_{1}+\mathbf{...}%
+\mathbf{a}_{d}s_{d}+\mathbf{e}_{d},  \label{mp1}
\end{equation}%
which is a 3D vector of length $R_{d}$ (the scalar product here is taken in $\mathbb{R}^{3}$), 
\begin{equation}
R_{d}\left( s_{1},...,s_{d}\right) =\left| \mathbf{R}_{d}\left(
s_{1},...,s_{d}\right) \right| =\sqrt{\mathbf{R}_{d}\left(
s_{1},...,s_{d}\right) \cdot \mathbf{R}_{d}\left( s_{1},...,s_{d}\right) }.
\label{mp2}
\end{equation}%
Vector $\mathbf{e}_{d}$ can be decomposed into components parallel ($\mathbf{%
e}_{d\mathbf{\Vert }}$) and perpendicular ($\mathbf{e}_{d\mathbf{\perp }}$)
to the vector space defined by vectors $\mathbf{a}_{1},\mathbf{...},\mathbf{a%
}_{d}$, 
\begin{equation}
\mathbf{e}_{d}=\mathbf{e}_{d\mathbf{\Vert }}+\mathbf{e}_{d\mathbf{\perp }%
},\quad \mathbf{e}_{d\mathbf{\Vert }}\cdot \mathbf{e}_{d\mathbf{\perp }%
}=0,\quad \mathbf{e}_{d\mathbf{\Vert }}\in span\left( \mathbf{a}_{1},\mathbf{%
...},\mathbf{a}_{d}\right) .  \label{mp3}
\end{equation}%
The latter statement means,%
\begin{equation}
\mathbf{e}_{d\mathbf{\Vert }}=s_{10}\mathbf{a}_{1}+...+s_{d0}\mathbf{a}_{d},
\label{mp4}
\end{equation}%
where $s_{10},...,s_{d0}$ are some coefficients. Note that $\mathbf{e}_{d\mathbf{\perp }}=\mathbf{0}$ for $d\geqslant 3$ 
if $\mathbf{a}_{1}\mathbf{,...,a}_{d}$ are not coplanar. Depending on $d$ and $\mathbf{a}_{1}\mathbf{,...,a}_{d}$, 
the $s_{10},...,s_{d0}$ may not be unique (Appendix gives details how to construct projections of vectors for a given 
$\left\{ \mathbf{a}_{1},\mathbf{...},\mathbf{a}_{d}\right\} $). 
Both these facts do not prevent representation of $\mathbf{R}$ from (\ref{mp1}) in the form
\begin{eqnarray}
\mathbf{R}_{d}\left( s_{1},...,s_{d}\right) &=&\mathbf{P}_{d}\left(
s_{1},...,s_{d}\right) +\mathbf{e}_{d\bot },\quad \mathbf{P}_{d}\cdot 
\mathbf{e}_{d\bot }=0,\quad i=1,...,d,  \label{mp5} \\
\mathbf{P}_{d}\left( s_{1},...,s_{d}\right) &=&\mathbf{a}_{1}\left(
s_{1}+s_{10}\right) +\mathbf{...}+\mathbf{a}_{d}\left( s_{d}+s_{d0}\right) ,
\notag \\
R_{d}^{2} &=&P_{d}^{2}+h_{d}^{2},\quad P_{d}=\left| \mathbf{P}_{d}\left(
s_{1},...,s_{d}\right) \right| ,\quad h_{d}=\left| \mathbf{e}_{d\bot
}\right| .  \notag
\end{eqnarray}
Consider the integral $I_{d}$   of an arbitrary integrable function $G_{d}$ over a d-dimensional polyhedron $D_{d}$ with $\left(d-1\right)$-dimensional boundary $B_{d-1}$: 
\begin{equation}
I_{d}=\int_{D_{d}}G_{d}\left( R_{d}\left( s_{1},...,s_{d}\right) \right)
dD_{d}=\int_{D_{d}}G_{d}\left( R_{d}\left( s_{1},...,s_{d}\right) \right)
ds_{1}...ds_{d}.  \label{mp6}
\end{equation}
 We reduce this to a sum of ($d-1$) dimensional integrals using the divergence theorem,%
\begin{equation}
I_{d}=\sum_{j}I_{\left( d-1\right) j},\quad I_{(d-1)j}=\int_{B_{(d-1)j}}%
\mathbf{N}_{dj}\odot \mathbf{V}_{d}dB_{d-1},\quad G_{d}\left( R_{d}\right)
=\nabla _{d}\odot \mathbf{V}_{d}.  \label{mp7}
\end{equation}%
Here $\odot $ is used for the dot product in $d$ dimensional space defined by $\mathbf{i}_{1},...,\mathbf{i}_{d}$ 
(to avoid confusion with the dot product in $\mathbb{R}^3$ defined by vectors $\mathbf{i}_{x},\mathbf{i}_{y},\mathbf{i}_{z}$). 
We assume that boundary $B_{d-1}$ consists of $\left( d-1\right) $ dimensional ``faces'' $B_{\left( d-1\right)j}$, which have external $d$-dimensional normal $\mathbf{N}_{dj}$, and $%
\mathbf{V}_{d}$ is a $d$-dimensional vector field which should be
constructed.
Certainly, $\mathbf{V}_{d}$ \ is not unique. The following procedure is applicable for any dimension and can be repeated recursively as the dimensionality reduces. 
It is based on the Euler theorem on homogeneous functions~\cite{KornKorn}. Since $P_{d}^{2}$ is a homogeneous function of degree 2 of $s_{1}+s_{10},...,s_{d}+s_{d0}$, we have  
\begin{equation}
P_{d}^{2}=\frac{1}{2}\left[ \left( s_{1}+s_{10}\right) \frac{\partial
P_{d}^{2}}{\partial s_{1}}+...+\left( s_{d}+s_{d0}\right) \frac{\partial
P_{d}^{2}}{\partial s_{d}}\right] .  \label{mp8}
\end{equation}%
This relation can be checked directly. Thus, if $\mathbf{V}_{d}$ is sought in the form%
\begin{equation}
\mathbf{V}_{d}=\left[ \mathbf{i}_{1}\left( s_{1}+s_{10}\right) +...+\mathbf{i%
}_{d}\left( s_{d}+s_{d0}\right) \right] F_{d}\left( P_{d}\right) ,
\label{mp9}
\end{equation}%
where we call $F_{d}\left( P_{d}\right) $ as the primitive
boundary function (PBF), we have 
\begin{equation}\begin{aligned}
\nabla _{d}\odot \mathbf{V}_{d} =&\left[ \mathbf{i}_{1}\left(
s_{1}+s_{10}\right) +...+\mathbf{i}_{d}\left( s_{d}+s_{d0}\right) \right]
\odot \nabla _{d}F_{d}\left( P_{d}\right) \\&+ F_{d}\left( P_{d}\right) \nabla
_{d}\odot \left[ \mathbf{i}_{1}\left( s_{1}+s_{10}\right) +...+\mathbf{i}%
_{d}\left( s_{d}+s_{d0}\right) \right]   \label{mp10} \\
=&\left[ \left( s_{1}+s_{10}\right) \frac{\partial P_{d}}{\partial s_{1}}%
+...+\left( s_{d}+s_{d0}\right) \frac{\partial P_{d}}{\partial s_{d}}\right] 
\frac{dF_{d}\left( P_{d}\right) }{dP_{d}}+F_{d}\left( P_{d}\right) d \\
=&P_{d}\frac{dF_{d}\left( P_{d}\right) }{dP_{d}}+F_{d}\left( P_{d}\right) d.
\end{aligned}\end{equation}%
Due to the last of (\ref{mp7}) and (\ref{mp5}) between $R_{d}$ and $P_{d}$,  $F_{d}$ satisfies
\begin{equation}
P_{d}\frac{dF_{d}\left( P_{d}\right) }{dP_{d}}+F_{d}\left( P_{d}\right)
d=G_{d}\left( \sqrt{P_{d}^{2}+h_{d}^{2}}\right) .  \label{mp11}
\end{equation}
Multiplying both sides of this differential equation by $P_{d}^{d-1}$ the left hand side can be expressed as the derivative of function $P_{d}^{d}F_{d}\left( P_{d}\right) $. 
Integrating, we get
\begin{equation}
F_{d}\left( P_{d}\right) =\frac{1}{P_{d}^{d}}\int P_{d}^{d-1}G_{d}\left( 
\sqrt{P_{d}^{2}+h_{d}^{2}}\right) dP_{d}+C_{d}.  \label{mp12}
\end{equation}%
The constant $C_{d}$ can be selected to avoid singularity of $F_{d}$ at $h_{d}\neq 0$ and $P_{d}=0$:%
\begin{equation}
F_{d}\left( P_{d}\right) =\frac{1}{P_{d}^{d}}\int_{0}^{P_{d}}P^{d-1}G_{d}%
\left( \sqrt{P^{2}+h_{d}^{2}}\right) dP,  \label{mp14}
\end{equation}%
which fixes the constant of the  PBF. 
A special case $h_{d}=0$ ($P_{d}=R_{d}$), $G_{d}(R_{d})=R_{d}^{-\alpha }$ is of interest.  We have the solution of (\ref{mp11}) as
\begin{equation}\begin{aligned}
F_{d}\left( P_{d}\right) =\left\{ 
\begin{array}{c}
\frac{1}{d-\alpha }P_{d}^{-\alpha }=\frac{1}{d-\alpha }G_{d}(P_{d})=\frac{1}{%
d-\alpha }G_{d}(R_{d}),\quad d\neq \alpha.  \\ 
P_{d}^{-\alpha }\ln P_{d}=R_{d}^{-\alpha }\ln R_{d},\quad d=\alpha .%
\end{array}%
\right.  \label{mp15}
\end{aligned}\end{equation}
The $\left( d-1\right) $-dimensional ``faces'' $B_{\left( d-1\right) j}$ considered are simple:  
points ($d=1$), line segments ($d=2$), rectangular triangles or rectangles ($d=3$), and right triangular prisms ($d=4$). 
In $d$-dimensional space these objects may correspond only to hyperplanes $s_{i}=0$, $s_{i}=1$ ($\mathbf{N}_{d_{j}}=\mp \mathbf{i}_{i}$) and 
$s_{i}+s_{k}=1$ ($\mathbf{N}_{d_{j}}=\left( \mathbf{i}_{i}+\mathbf{i}_{k}\right) /\sqrt{2}$), $i,k=1,...,d$ (for pairs $i,k$ describing parametrization of the same triangle). 
Using (\ref{mp9}) we obtain for (\ref{mp7})%
\begin{equation}\begin{aligned}
I_{(d-1)j}
=&\left\{ 
\begin{array}{c}\!\!\!\!
-s_{i0}\int_{B_{(d-1)j}}F_{d}\left( \left. P_{d}\right| _{s_{i}=0}\right)
ds_{1}..ds_{i-1}ds_{i+1}..ds_{d},\quad s_{i}=0, \\ 
\left( 1+s_{i0}\right) \int_{B_{(d-1)j}}F_{d}\left( \left. P_{d}\right|
_{s_{i}=1}\right) ds_{1}..ds_{i-1}ds_{i+1}..ds_{d},\quad s_{i}=1, \\ 
\!\!\!\! \left( 1\!+\!s_{i0}\!+\!s_{k0}\right) \int_{B_{(d-1)j}} \!\! F_{d}\left( \left.
P_{d}\right| _{s_{i}+s_{k}=1}\right) ds_{1}..ds_{i-1}ds_{i+1}..ds_{d},\quad \!\!\!\!
s_{i}+s_{k}=1.%
\end{array}%
\right.  \label{mp16}
\end{aligned}\end{equation}%
Here we used integration over $s_i$ instead of the triangle hypotenuse $s_{i}+s_{k}=1$. 
This can be done for any of the integrals in $\left( d-1\right) $-dimensional space with coordinates 
\begin{equation}
s_{1}^{\prime }=s_{1},...,s_{i-1}^{\prime }=s_{i-1},s_{i}^{\prime
}=s_{i+1},...,s_{d-1}^{\prime }=s_{d}.  \label{mp16.1}
\end{equation}%
The $P_{d}$ in the integrals (\ref{mp16}) can be denoted as $R_{d-1}$ and $F_{d}$ as $G_{d-1}$. 
The latter can be brought to standard form (\ref{mp1}), e.g., for the last integral in (\ref{mp16}) we have
\begin{gather}
\left. P_{d}\right| _{s_{i}+s_{k}=1} \! = \! 
\left| \sum_{l\neq i,k}\mathbf{a}_{l}s_{l}+\mathbf{a}_{i}\left( 1-s_{k}\right) +\mathbf{a}_{k}s_{k}+\mathbf{e}_{d\mathbf{\Vert }}\right| \!=\!\left| \sum_{l=1}^{d-1}\mathbf{a}_{l}^{\prime}s_{l}^{\prime }+\mathbf{e}_{d-1}\right|\! =\!R_{d-1} \! \label{mp16.2} \\
\mathbf{a}_{l}^{\prime } = \left\{ 
\begin{array}{c}
\mathbf{a}_{l},\text{ }l<i,\quad l\neq k \\ 
\mathbf{a}_{l+1},\text{ }l>i,\quad l\neq k \\ 
\mathbf{a}_{k}-\mathbf{a}_{i},\quad l=k%
\end{array}%
\right. ,\quad \mathbf{e}_{d-1}=\mathbf{e}_{d\mathbf{\Vert }}+\mathbf{a}_{i}.  
\end{gather}%
\begin{figure}[tbh]
	\begin{center}
		\includegraphics[width=0.9\textwidth]{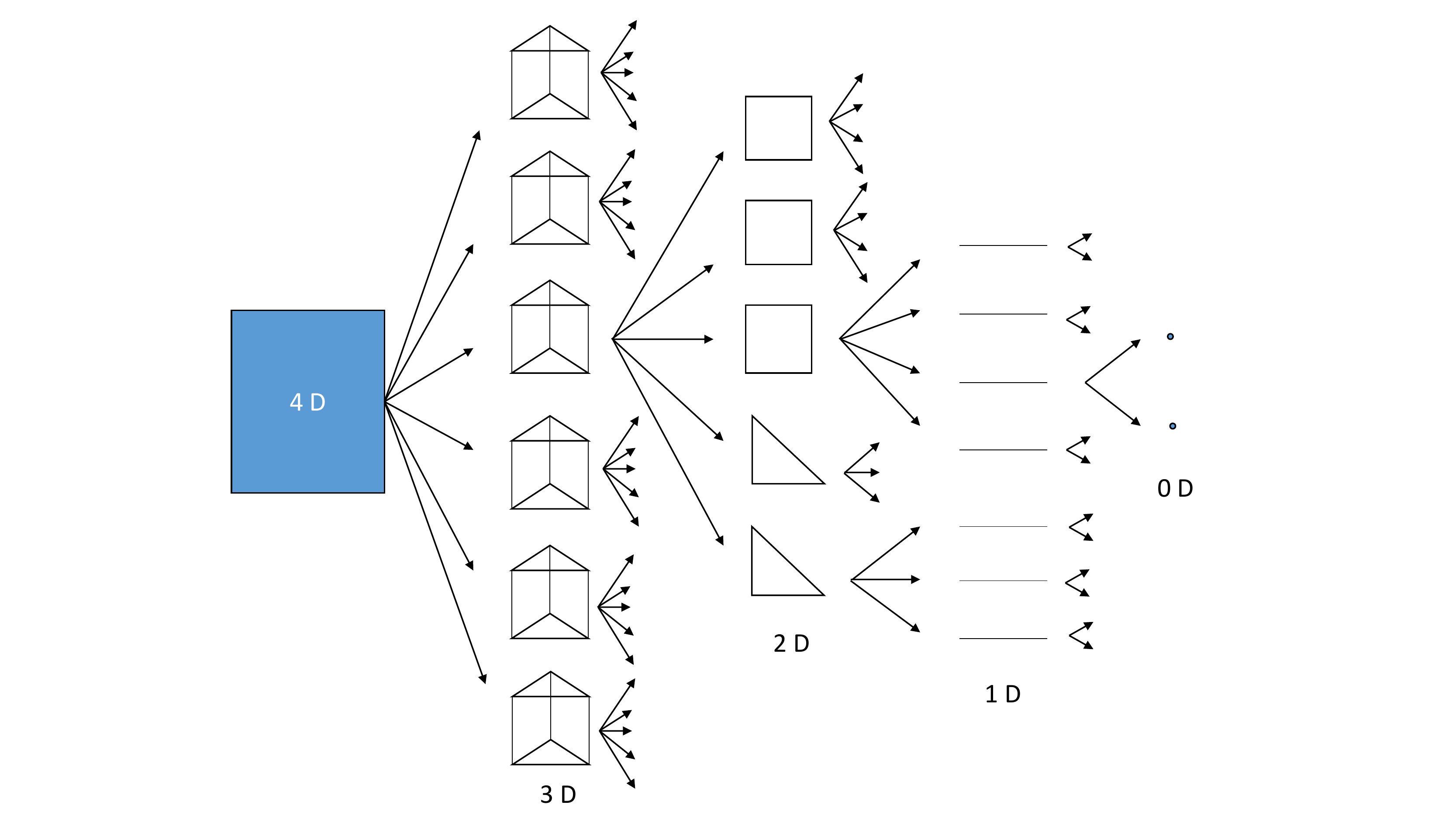}
	\end{center}
     \setlength{\abovecaptionskip}{-3pt}
	\caption{A scheme of recursive reduction of 4D integrals to sums of lower
		dimensional integrals. }
	\label{Fig1}
\end{figure}
Fig. \ref{Fig1} shows the scheme, where the Gauss' divergence theorem is applied at each dimension reduction step. 
Note that while $h_{4}$ is common to all domains,  $h_{d}$ for $d<4$ (as well as $\mathbf{a}_{i}$ and $\mathbf{e}_{d}$) 
are specific to each domain. ``Child'' domains of dimension $d$ inherit the values of parameters $h_{d+1},h_{d+2},...$ 
from ``parents'' and add domain specific $h_{d}$ to the list, which they pass to ``children.'' 
We  use  permutation of variables $s_{1},...,s_{d}$ to reduce the number of integrals over faces $B_{\left( d-1\right) j}$ to a minimum.  
Integrand transformations on dimensionality reduction can be described as%
\begin{equation}
G_{d-1}\left( R_{d-1}\right) =F_{d}\left( P_{d}\right) ,\quad
R_{d-1}^{2}=P_{d}^{2}=R_{d}^{2}-h_{d}^{2},\quad d=4,3,2,1.  \label{mp17}
\end{equation}%
With (\ref{mp5}) and (\ref{mp14}) a recursive process is created, and the final PBF is 
\begin{gather}
\label{mp18} F_{1}\left( P_{1};h_{1},h_{2},h_{3},h_{4}\right)= \\
\notag \frac{1}{P_{1}}%
\int_{0}^{P_{1}} \!\!\!\!\!\! \frac{1}{p_{1}^{2}+h_{1}^{2}}
\int_{0}^{\sqrt{p_{1}^{2}+h_{1}^{2}}} \!\!\!\!\!\!\!\!\!\!\!\! \frac{p_{2}}{\left( p_{2}^{2}+h_{2}^{2}\right) ^{\frac{3}{2}}}
\int_{0}^{\sqrt{p_{2}^{2}+h_{2}^{2}}} \!\!\!\!\!\!\!\!\!\! \frac{p_{3}^{2}}{\left(
p_{3}^{2}+h_{3}^{2}\right) ^{2}}
\int_{0}^{\sqrt{p_{3}^{2}+h_{3}^{2}}} \!\!\!\!
p_{4}^{3}G_{4}\left( \sqrt{p_{4}^{2}+h_{4}^{2}}\right)dp_{4}...
dp_{1}.  
\end{gather}%
All integrals for the free space Green's function of the
Laplace equation can be calculated analytically. At the end of the reduction process the 1D domains are nothing but unit
segments $s_{1}\in \left[ 0,1\right] $. From (\ref{mp7}) and (%
\ref{mp16}) we have%
\begin{gather}
\notag I_{1}\left( h_{1},h_{2},h_{3},h_{4}\right)  =\left( 1+s_{10}\right) F_{1}\left( P_{1}\left( 1\right) ;h_{1},h_{2},h_{3},h_{4}\right) 
-s_{10}F_{1}\left( P_{1}\left( 0\right) ;h_{1},h_{2},h_{3},h_{4}\right) ,  \\
P_{1}\left( s_{1}\right)  =\left| s_{1}+s_{10}\right| \left| \mathbf{a}%
_{1}\right| . \label{mp19} 
\end{gather}%
The recursive process and analytical expressions for  (\ref{mp14}) or (\ref{mp18}) are described next.

\subsection{Primitive boundary functions}

We provide analytical expressions for PBF $F_{d}\left(
P_{d};h_{d},...,h_{4}\right) $. Besides $P_{d}$ these functions
depend on $h_{d},...,h_{4}$ only. 
Some of the parameters are zero, which leads to a substantial simplification of the expressions. 
We first prove several simple lemmas and specify all possible cases, and then provide expressions for PBFs for each case.
\subsubsection{Parameters $h_{d}$}
\begin{lemma}
In the set of parameters $\left\{ h_{1},...,h_{4}\right\} $ at least two
parameters are zero.
\end{lemma}
\begin{proof}
The method the projections are constructed (see (\ref{mp1}) and (\ref{mp3})) shows
\begin{equation}
\mathbf{R}_{4}=\mathbf{e}_{4\mathbf{\perp }}+\mathbf{R}_{3}=\mathbf{e}_{4%
\mathbf{\perp }}+\mathbf{e}_{3\bot }+\mathbf{R}_{2}=...=\mathbf{e}_{4\mathbf{%
\perp }}+\mathbf{e}_{3\bot }+\mathbf{e}_{2\perp }+\mathbf{e}_{1\bot }+%
\mathbf{P}_{1}, \label{b0}
\end{equation}%
with $h_{d}=\left\| \mathbf{e}_{d\bot }\right\| ,\quad
d=1,...,4$.
By construction all $\mathbf{e}_{d\bot }$ are perpendicular to each other 
and also to $\mathbf{P}_{1}$. Since these five vectors lie in 3D space, two of them should be zero. 
However, $P_{1}=\left\| \mathbf{P}_{1}\right\| \neq 0$ (as the final integration is taken 
over some edge - a segment of finite length), so, at least two of vectors $\mathbf{e}_{d\bot }$ are zero.
\end{proof}

\begin{lemma}
In the set  $\left\{ h_{1},h_{2},h_{3}\right\} $ at least one
is nonzero. So, $h_{1}+h_{2}+h_{3}\neq 0.$
\end{lemma}

\begin{proof}
Independent of $h_{4}$, $\mathbf{R}_{3}$ has at least two non-zero
components belonging to $span\left( \mathbf{X}_{s},\mathbf{X}_{t}\right) $
or to $span\left( \mathbf{Y}_{u},\mathbf{Y}_{v}\right) $. As $\mathbf{R}_{3}$
can be expanded over 4 vectors and $\mathbf{P}_{1}$ in (\ref{b0}) is not
zero, at least one of $\mathbf{e}_{d\bot }$, $d=1,...,3$ is not
zero. 
\end{proof}

\begin{remark}
As follows from these two lemmas either two or three parameters $\left\{
h_{1},...,h_{4}\right\} $ are zero.
\end{remark}

\begin{lemma}
$h_{3}h_{4}=0$.
\end{lemma}

\begin{proof}
Indeed, if $h_{4}\neq 0$ then the planes $P_{x}$ and $P_{y}$ containing
triangles $S_{x}$ and $S_{y}$ are parallel. In this case, $\mathbf{e}_{3}=$ $%
\mathbf{e}_{4\mathbf{\Vert }}\in span\left( \mathbf{X}_{s},\mathbf{X}%
_{t}\right) =span\left( \mathbf{Y}_{u},\mathbf{Y}_{v}\right) $. Since any
three out of four vectors $\mathbf{X}_{s},\mathbf{X}_{t}$,$\mathbf{Y}_{u},%
\mathbf{Y}_{v}$ generate exactly $span\left( \mathbf{X}_{s},\mathbf{X}%
_{t}\right) $, then $\mathbf{e}_{3||}$ belongs to the same subspace as $%
\mathbf{e}_{3}$ and $\mathbf{e}_{3\bot }=\mathbf{0}$.
\end{proof}

\begin{lemma}
$h_{1}+h_{2}\neq 0.$
\end{lemma}

\begin{proof}
If $h_{3}=0$ then Lemma (2) provides $h_{1}+h_{2}\neq 0$. If $h_{3}\neq 0$
then according to Lemma (3) we have $h_{4}=0$. The latter shows that
planes $P_{x}$ and $P_{y}$ containing triangles $S_{x}$ and $S_{y}$ are not
parallel and $span\left( \mathbf{X}_{s},\mathbf{X}_{t},\mathbf{Y}_{u},%
\mathbf{Y}_{v}\right) $ $=\mathbb{R}^{3}$. So, in representation (\ref{b0})
of $\mathbf{R}_{4}$ at least three out of five $\left\{ \mathbf{e}_{4\bot },%
\mathbf{e}_{3\bot },\mathbf{e}_{2\bot },\mathbf{e}_{1\bot },\mathbf{P}%
_{1}\right\} $ are nonzero. As $\mathbf{e}_{4\bot }=\mathbf{0},\mathbf{e}%
_{3\bot }\neq \mathbf{0}$, and $\mathbf{P}_{1}\neq \mathbf{0}$ one of the
vectors $\left\{ \mathbf{e}_{2\bot },\mathbf{e}_{1\bot }\right\} $ should be
nonzero. 
\end{proof}
As a result only the following cases of parameters $\left\{ h_{1},...,h_{4}\right\} $ are possible: 
\begin{equation}\begin{aligned}
\#1\quad h_{4} &= h_{3}=h_{2}=0,\quad h_{1}\neq 0,\quad  \label{b0.1} \quad \quad \#2\quad h_{4} = h_{3}=h_{1}=0,\quad h_{2}\neq 0,  \\
\quad \#3\quad h_{4} &= h_{3}=0,\quad h_{2}\neq 0,\quad h_{1}\neq 0,  \quad  \#4\quad h_{4} = h_{2}=0,\quad h_{3}\neq 0,\quad h_{1}\neq 0,\\ 
\#5\quad h_{4} &= h_{1}=0,\quad h_{3}\neq 0,\quad h_{2}\neq 0,  \quad   \#6\quad h_{3} = h_{2}=0,\quad h_{4}\neq 0,\quad h_{1}\neq 0,  \\
\#7\quad h_{3} &= h_{1}=0,\quad h_{4}\neq 0,\quad h_{2}\neq 0.\quad \quad 
\end{aligned}\end{equation}%

\begin{figure}[tbh]
 \vspace*{-0.1in}
	\begin{center}
		\includegraphics[width=0.9\textwidth]{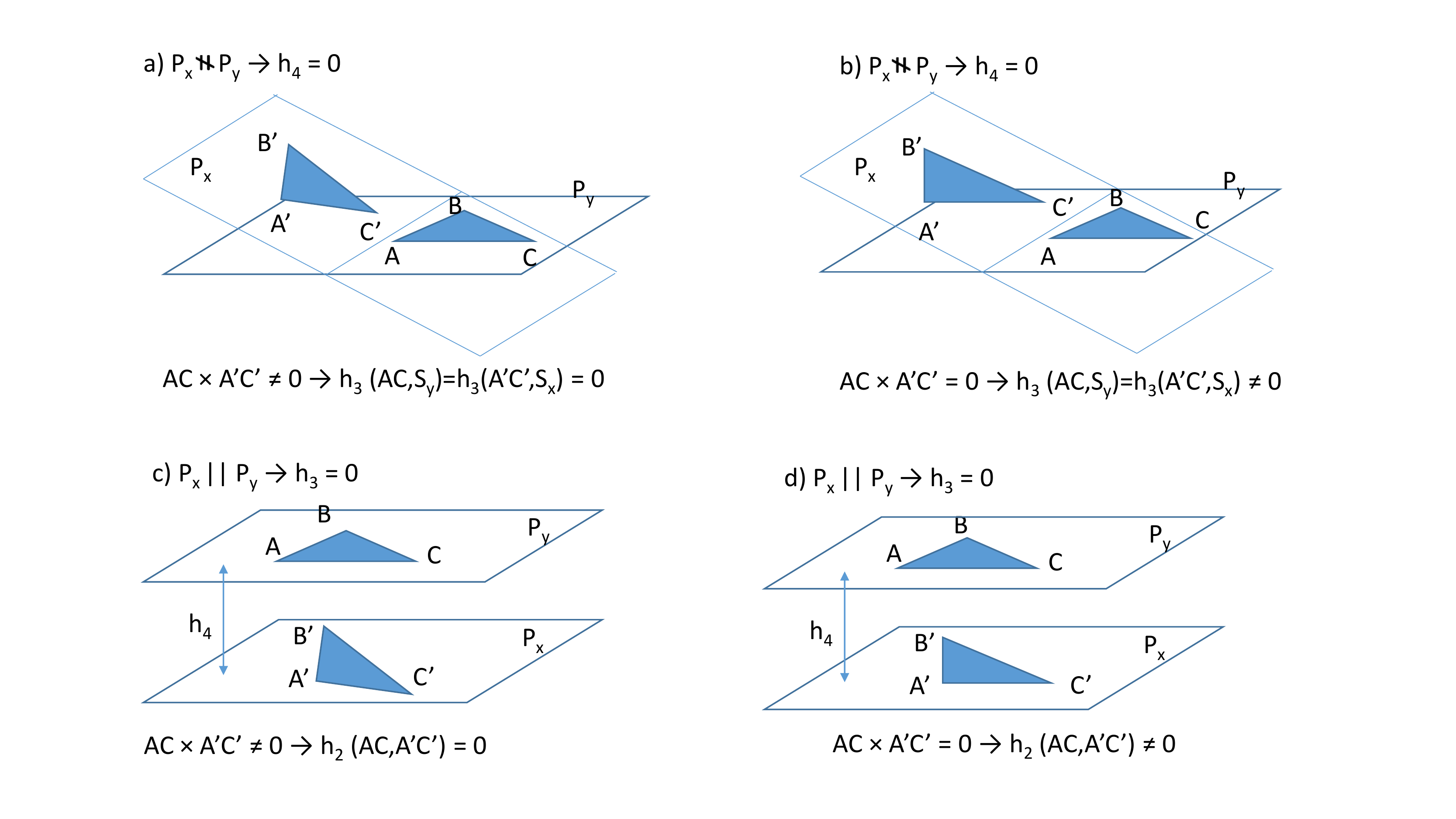}
	\end{center}
     \setlength{\abovecaptionskip}{-10pt}
     \setlength{\belowcaptionskip}{-20pt}
	\caption{ Illustration of typical cases determined by the location of the source ($S_{x}$) and receiver ($S_{y}$) triangles residing in planes $P_{x}$ and $P_{y}$, respectively.}
	\label{Fig2}
\end{figure}

Some of the typical situations are illustrated in Fig. \ref{Fig2}.

\subsubsection{PBF for 4D integral}

For the single layer potential we have%
\begin{equation}
G_{4}\left( R_{4}\right) =R_{4}^{-1}.  \label{k1}
\end{equation}%
General expression for the PBF for $d=4$ is (see (\ref{mp14}))
\begin{equation}
F_{4}\left( P_{4};h_{4}\right) =\frac{1}{3P_{4}^{4}}\left[ \left(
P_{4}^{2}-2h_{4}^{2}\right) R_{4}+2h_{4}^{3}\right] ,\quad R_{4}=\sqrt{%
P_{4}^{2}+h_{4}^{2}}.  \label{k2}
\end{equation}%
We do not provide derivation of formulae like this, as all
obtained expressions can be reduced to tabulated integrals of
elementary functions. A simple way
to check the validity of the expressions is to differentiate them and check
that (\ref{mp11}) is satisfied (or differentiate $%
P_{d}^{d}F_{d}$ and check that the result is $G_{d}\left( R_{d}\right) $).
An important case is $h_{4}=0.$ It can be obtained by
setting $h_{4}=0$ in (\ref{k2}), or by using (\ref{mp15}),%
\begin{equation}
F_{4}\left( P_{4};0\right) =\frac{1}{3P_{4}}=\frac{1}{3}G_{4}(R_{4}),\quad
R_{4}=P_{4}.  \label{k3}
\end{equation}%
So (\ref{k3}) can be used for cases 1-5, while (\ref{k2}) for cases 6-7
in (\ref{b0.1}).

\subsubsection{PBF for 3D integrals}

Here, at least one of the parameters $h_{3}$ and $h_{4}$ is zero. Using (\ref%
{mp14}) and (\ref{mp15}) for $G_{3}\left( R_{3}\right) =F_{4}\left(
P_{4}\right) $ and $d=3$, we obtain for the cases listed in (\ref{b0.1})%
\begin{equation}\begin{aligned}
\#1,\#2,\#3 :\quad F_{3}\left( P_{3};0,0\right) &=\frac{1}{6P_{3}},
\label{k4} \\
\#4,\#5 :\quad F_{3}\left( P_{3};h_{3},0\right) &=\frac{1}{6P_{3}^{3}}%
\left( P_{3}R_{3}-h_{3}^{2}\ln \frac{P_{3}+R_{3}}{h_{3}}\right) ,\quad R_{3}=%
\sqrt{P_{3}^{2}+h_{3}^{2}}, \\
\#6,\#7 :\quad F_{3}\left( P_{3};0,h_{4}\right) &=\frac{1}{6P_{3}^{3}}%
\left( P_{3}R_{4}-3h_{4}^{2}\ln \frac{P_{3}+R_{4}}{h_{4}}+4h_{4}^{2}\frac{%
R_{4}-h_{4}}{P_{3}}\right) ,\\
R_{4}&=\sqrt{P_{3}^{2}+h_{4}^{2}}. 
\end{aligned}\end{equation}%

\subsubsection{PBF for 2D integrals}

Now in (\ref{k4}) we substitute $G_{2}\left( R_{2}\right) =F_{3}\left(
P_{3}\right) $ and use (\ref{mp14}) and (\ref{mp15}) at $d=2$ to obtain the
following expressions for the cases listed in (\ref{b0.1}) 
 
\begin{equation}\begin{aligned}
&\#1:\quad F_{2}\left( P_{2};0,0,0\right) =\frac{1}{6P_{2}},  \label{k5} \quad
\#2,\#3:\quad F_{2}\left( P_{2};h_{2},0,0\right) =\frac{R_{2}-h_{2}}{%
6P_{2}^{2}},  \\
&\#4:\quad F_{2}\left( P_{2};0,h_{3},0\right) =\frac{1}{6P_{2}^{2}}\left[
R_{3}-2h_{3}+\frac{h_{3}^{2}}{P_{2}}\ln \frac{P_{2}+R_{3}}{h_{3}}\right]
,\quad  \\
&\#5:\quad F_{2}\left( P_{2};h_{2},h_{3},0\right) =\frac{1}{6P_{2}^{2}}\left[
R_{3}-C_{5}+\frac{h_{3}^{2}}{R_{2}}\ln \frac{R_{2}+R_{3}}{h_{3}}\right] , 
\\
&\#6:\quad F_{2}\left( P_{2};0,0,h_{4}\right) =\frac{1}{6P_{2}^{2}}\left[
R_{4}-3h_{4}+h_{4}^{2}\left( \frac{3}{P_{2}}\ln \frac{P_{2}+R_{4}}{h_{4}}-2%
\frac{R_{4}-h_{4}}{P_{2}^{2}}\right) \right] ,\quad   \\
&\#7:\quad F_{2}\left( P_{2};h_{2},0,h_{4}\right) =\frac{1}{6P_{2}^{2}}\left[
R_{4}-C_{7}+h_{4}^{2}\left( \frac{3}{R_{2}}\ln \frac{R_{2}+R_{4}}{h_{4}}-2%
\frac{R_{4}-h_{4}}{R_{2}^{2}}\right) \right] ,  \\
&h=\sqrt{h_{2}^{2}+h_{3}^{2}+h_{4}^{2}},\quad  R_{3}=\sqrt{P_{2}^{2}+h_{2}^{2}+h_{3}^{2}},\quad R_{4}=\sqrt{P_{2}^{2}+h_{2}^{2}+h_{3}^{2}+h_{4}^{2}}, \\
&R_{2}=\sqrt{P_{2}^{2}+h_{2}^{2}},\quad C_{5}=h+\frac{h_{3}^{2}}{h_{2}}\ln \frac{h_{2}+h}{h_{3}},\quad \!\! \!
C_{7}=h+h_{4}^{2}\left( \frac{3}{h_{2}}\ln \frac{h_{2}+h}{h_{4}}-2\frac{%
h-h_{4}}{h_{2}^{2}}\right) .   
\end{aligned}\end{equation}%

\subsubsection{PBF for 1D integrals}

Similarly, by setting $G_{1}\left( R_{1}\right) =F_{2}\left( P_{2}\right) $
and using (\ref{mp14}) at $d=1$ we obtain 
\begin{equation}\begin{aligned}
&\#1 :\quad F_{1}\left( P_{1};h_{1},0,0,0\right) =\frac{1}{6}\Phi
_{1}\left( P_{1}\right) ,  \label{k6} \\&
\#2 :\quad F_{1}\left( P_{1};0,h_{2},0,0\right) =\frac{1}{6}\left( \Phi
_{1}\left( P_{1}\right) -\frac{1}{R_{2}+h_{2}}\right) ,  \\
&\#3 :\quad F_{1}\left( P_{1};h_{1},h_{2},0,0\right) =\frac{1}{6}\left[
\Phi _{1}\left( P_{1}\right) -\frac{h_{2}}{h_{1}}\Phi _{2}\left(
P_{1};h_{1}\right) \right] , \\
&\#4 :\quad F_{1}\left( P_{1};h_{1},0,h_{3},0\right) =\frac{1}{6}%
\left[ \left( 1-\frac{h_{3}^{2}}{h_{1}^{2}}\right) \Phi _{1}\left(
P_{1}\right) -\frac{2h_{3}}{h_{1}}\Phi _{2}\left( P_{1};h_{1}\right) +\frac{%
h_{3}^{2}}{h_{1}^{2}}\Phi _{3}\left( P_{1}\right) \right] ,\quad   \\
&\#5 :\quad F_{1}\left( P_{1};0,h_{2},h_{3},0\right) =\frac{1}{6}\left[ 
\frac{h^{2}}{h_{2}^{2}}\Phi _{1}\left( P_{1}\right) -\frac{1}{R_{4}+h}-\Phi
_{4}\left( P_{1};h_{3}\right) \right] ,\quad   \\
&\#6 :\quad F_{1}\left( P_{1};h_{1},0,0,h_{4}\right) = \\ & \quad \frac{1}{6}%
\left[ \left( 1-\frac{3h_{4}^{2}}{h_{1}^{2}}\right) \Phi _{1}\left(
P_{1}\right) -\left( 3-\frac{h_{4}^{2}}{h_{1}^{2}}\right) \frac{h_{4}}{h_{1}}%
\Phi _{2}\left( P_{1};h_{1}\right) +3\frac{h_{4}^{2}}{h_{1}^{2}}\Phi
_{3}\left( P_{1}\right) -\frac{h_{4}^{2}}{h_{1}^{2}\left( R_{4}+h_{4}\right) 
}\right] ,   \\
&\#7 :\quad F_{1}\left( P_{1};0,h_{2},0,h_{4}\right) \\&\quad\quad\quad=\frac{1}{6}\left[
\left( 1+\frac{3h_{4}^{2}}{h_{2}^{2}}\right) \Phi _{1}\left( P_{1}\right) -%
\frac{2h_{4}^{3}}{h_{2}^{3}}\Phi _{2}\left( P_{1};h_{2}\right) -3\Phi
_{4}\left( P_{1};h_{4}\right) +\frac{2h_{4}^{2}-h_{2}^{2}}{h_{2}^{2}\left(
R_{4}+h\right) }\right] ,
\end{aligned}\end{equation}
\begin{equation}\begin{aligned}
&R_{4} =\sqrt{P_{1}^{2}+h^{2}},\quad R_{1}=\sqrt{P_{1}^{2}+h_{1}^{2}},\quad
R_{2}=\sqrt{P_{1}^{2}+h_{2}^{2}},\quad
h^{2}=h_{1}^{2}+h_{2}^{2}+h_{3}^{2}+h_{4}^{2},   \\
&\Phi _{1}\left( P_{1}\right) =\frac{1}{P_{1}}\ln \frac{P_{1}+R_{4}}{h}%
,\quad \Phi _{2}\left( P_{1};\eta \right) =\frac{1}{P_{1}}\arctan \frac{\eta
P_{1}}{h^{2}+R_{4}\sqrt{h^{2}-\eta ^{2}}},  \\
&\Phi _{3}\left( P_{1}\right) =\frac{1}{R_{1}}\ln \frac{R_{1}+R_{4}}{\sqrt{%
h^{2}-h_{1}^{2}}},\quad \Phi _{4}\left( P_{1};\eta \right) =\frac{\eta ^{2}}{%
P_{1}^{2}h_{2}}\left( \frac{R_{2}}{h_{2}}\ln \frac{R_{2}+R_{4}}{\eta }-\ln 
\frac{h_{2}+h}{\eta }\right) .  \notag
\end{aligned}\end{equation}

\subsection{Domains and parameters}

In this section we specify the domains of integration and parameters, which
enable determination of projections and parameters $h_{1},...,h_{4}$. Before
doing that we provide a proof of a lemma, which allows us obtain substantial
simplifications of analytical expressions.

\subsubsection{4D domain and 4D integral}

From the parametrization of triangles and (\ref{mp1}) and (\ref{mp2}) we
have for $d=4$%
\begin{equation}\begin{aligned}
s_{1} =& s,\quad s_{2}=t,\quad s_{3}=u,\quad s_{4}=v,  \label{b01} \\
\mathbf{a}_{1} =& \mathbf{X}_{s},\text{\quad }\mathbf{a}_{2}=\mathbf{X}%
_{t},\quad \mathbf{a}_{3}=-\mathbf{Y}_{u},\quad \mathbf{a}_{4}=-\mathbf{Y}%
_{v},\quad \mathbf{e}_{4}=\mathbf{X}_{0}-\mathbf{Y}_{0}.  
\end{aligned}\end{equation}%
We are interested in evaluation of integral%
\begin{gather}
U=I_{4}\left( \mathbf{a}_{1},\mathbf{a}_{2},\mathbf{a}_{3},\mathbf{a}_{4},%
\mathbf{e}_{4}\right)
=\int_{D_{4}}... \label{b02} \\ \notag =\int_{0}^{1}\int_{0}^{1-s_{1}}\int_{0}^{1}%
\int_{0}^{1-s_{3}}G_{4}\left( R_{4}\left( s_{1},s_{2},s_{3},s_{4};\mathbf{a}%
_{1},\mathbf{a}_{2},\mathbf{a}_{3},\mathbf{a}_{4},\mathbf{e}_{4}\right)
\right) ds_{4}ds_{3}ds_{2}ds_{1}  
\end{gather}%
The integration domain $D_{4}$ is the Cartesian product of standard triangles, $T_{x}\otimes $ $T_{y}$.

\subsubsection{4D domain boundaries and 3D integrals}

The boundaries of domain $D_{4}$ consist of 6 right prisms:
\begin{equation}\begin{aligned}
B_{31} &= T_{x}\otimes C_{y1},\quad B_{32}=T_{x}\otimes C_{y2},\quad
B_{33}=T_{x}\otimes C_{y3},  \label{b1} \\
B_{34} &= C_{x1}\otimes T_{y},\quad B_{35}=C_{x2}\otimes T_{y},\quad
B_{36}=C_{x3}\otimes T_{y}, 
\end{aligned}\end{equation}%
where $C_{xj}$ and $C_{yj}$, $j=1,2,3$ are the edges of triangles $T_{x}$
and $T_{y}$, respectively.  
(\ref{mp7}) and (\ref{mp16}) then provide%
\begin{equation}\begin{aligned}
U =&I_{4}\left( \mathbf{a}_{1},\mathbf{a}_{2},\mathbf{a}_{3},\mathbf{a}_{4},%
\mathbf{e}_{4}\right) =\sum_{j=1}^{6}I_{3j},  \label{b2} \\ 
I_{31} =&-s_{40}\int_{B_{31}}F_{4}\left( P_{4}\left(
s_{1},s_{2},s_{3},0\right) \right) ds_{1}ds_{2}ds_{3}=-s_{40}I_{3}\left( 
\mathbf{a}_{11},\mathbf{a}_{21},\mathbf{a}_{31},\mathbf{e}_{31}\right) , 
 \\
I_{32} =&\left( 1+s_{30}+s_{40}\right) \int_{B_{32}}F_{4}\left( P_{4}\left(
s_{1},s_{2},1-s_{4},s_{4}\right) \right) ds_{1}ds_{2}ds_{4} \\=& \left(
1+s_{30}+s_{40}\right) I_{3}\left( \mathbf{a}_{12},\mathbf{a}_{22},\mathbf{a}%
_{32},\mathbf{e}_{32}\right) ,   \\
I_{33} =&-s_{30}\int_{B_{33}}F_{4}\left( P_{4}\left(
s_{1},s_{2},0,s_{4}\right) \right) ds_{1}ds_{2}ds_{4}=-s_{30}I_{3}\left( 
\mathbf{a}_{13},\mathbf{a}_{23},\mathbf{a}_{33},\mathbf{e}_{33}\right) , 
\\
I_{34} =&-s_{20}\int_{B_{34}}F_{4}\left( P_{4}\left(
s_{1},0,s_{3},s_{4}\right) \right) ds_{1}ds_{3}ds_{4}=-s_{20}I_{3}\left( 
\mathbf{a}_{14},\mathbf{a}_{24},\mathbf{a}_{34},\mathbf{e}_{34}\right) , 
\\
I_{35} =& 
\left( 1+s_{10}+s_{20}\right) \int_{B_{35}}F_{4}\left( P_{4}\left(1-s_2,s_{2},s_{3},s_{4}\right) \right) ds_{2}ds_{3}ds_{4}
 \\=& \left(1+s_{10}+s_{20}\right) I_{3}\left( \mathbf{a}_{15},\mathbf{a}_{25},\mathbf{a}_{35},\mathbf{e}_{35}\right) ,  \\
I_{36} =&
-s_{10}\int_{B_{36}}F_{4}\left( P_{4}\left(0,s_{2},s_{3},s_{4}\right) \right) ds_{2}ds_{3}ds_{4}
=-s_{10}I_{3}\left( \mathbf{a}_{16},\mathbf{a}_{26},\mathbf{a}_{36},\mathbf{e}_{36}\right) .
\end{aligned}\end{equation}%
Here $D_{3}$ is a right prism in 3D space (the same as $B_{31}$), so 
\begin{equation}
I_{3}\left( \mathbf{a}_{1j},\mathbf{a}_{2j},\mathbf{a}_{3j},\mathbf{e}%
_{3j}\right)
=\!\!\int_{D_{3}}\!\!\!\!\!...=\!\!\int_{0}^{1}\!\!\!\int_{0}^{1-s_{1}}\!\!\!\!\int_{0}^{1}\!\!\!G_{3}\left(
R_{3}\left( s_{1},s_{2},s_{3};\mathbf{a}_{1j},\mathbf{a}_{2j},\mathbf{a}%
_{3j},\mathbf{e}_{3j}\right) \right) ds_{3}ds_{2}ds_{1}.  \label{b3}
\end{equation}%
Parameters of the 3D integral follow from the integration domain (to avoid
confusion we extend the subscript of parameters of these integrals to add
the index $j=1,...,6$ of the integral at the end),
\begin{equation}\begin{aligned}
\mathbf{a}_{11} &= \mathbf{a}_{1}\mathbf{,\quad a}_{21}=\mathbf{a}_{2}%
\mathbf{,\quad a}_{31}=\mathbf{a}_{3}\mathbf{,\quad a}_{12}=\mathbf{a}_{1}%
\mathbf{,\quad a}_{22}=\mathbf{a}_{2}\mathbf{,\quad a}_{32}=\mathbf{a}_{4}-%
\mathbf{a}_{3}\mathbf{,}  \label{b4} \\
\mathbf{a}_{13} &= \mathbf{a}_{1}\mathbf{,\quad a}_{23}=\mathbf{a}_{2}%
\mathbf{,\quad a}_{33}=\mathbf{a}_{4}\mathbf{,\quad a}_{14}=\mathbf{a}_{3}%
\mathbf{,\quad a}_{24}=\mathbf{a}_{4}\mathbf{,\quad a}_{34}=\mathbf{a}_{1}%
\mathbf{,}  \\
\mathbf{a}_{15} &= \mathbf{a}_{3}\mathbf{,\quad a}_{25}=\mathbf{a}_{4}%
\mathbf{,\quad a}_{35}=\mathbf{a}_{2}-\mathbf{a}_{1}\mathbf{,\quad a}_{16}=%
\mathbf{a}_{3}\mathbf{,\quad a}_{26}=\mathbf{a}_{4}\mathbf{,\quad a}_{36}=%
\mathbf{a}_{2}\mathbf{,}   \\ 
\mathbf{e}_{31} &= \mathbf{e}_{33}=\mathbf{e}_{34}=\mathbf{e}_{36}= \sum_{d=1}^{4} \mathbf{a}_{d}s_{d0},\quad 
\mathbf{e}_{32}=\mathbf{e}_{31}+\mathbf{a}_{3},\quad 
\mathbf{e}_{35}=\mathbf{e}_{31}+\mathbf{a}_{1},\quad 
\end{aligned}\end{equation}%
where the values for parameters $\mathbf{a}_{1},\mathbf{a}_{2},\mathbf{a}%
_{3},\mathbf{a}_{4}$ of the 4D integrals are provided by (\ref{b01}) and
coordinates of projections $s_{10},...,s_{40}$ and $h_{4}$ can be found from
these parameters and $\mathbf{e}_{4}$ (see Appendix).

\subsubsection{3D domain boundaries and 2D integrals}

The prism $D_{3}$ has 3 rectangular faces and 2 triangular faces.
Rectangles are located in planes $s_{1}+s_{2}=1$, $s_{1}=0$, \ and $s_{2}=0$%
, and triangles are located in planes $s_{3}=1$ and $s_{3}=0$ (see (\ref%
{b3})). So, we have%
\begin{equation}\begin{aligned}
&I_{3}\left( \mathbf{a}_{1},\mathbf{a}_{2},\mathbf{a}_{3},\mathbf{e}%
_{3}\right) =\sum_{j=1}^{5}I_{2j},  \label{b5} \\ 
&I_{21}=\left( 1+s_{10}+s_{20}\right) \int_{0}^{1}\int_{0}^{1}F_{3}\left(
P_{3}\left( s_{1},1-s_{1},s_{3}\right) \right) ds_{1}ds_{3}  = \left(
1+s_{10}+s_{20}\right) I_{2}^{(s)}\left( \mathbf{a}_{11},\mathbf{a}_{21},%
\mathbf{e}_{21}\right) ,  \\
&I_{22}=-s_{10}\int_{0}^{1}\int_{0}^{1}F_{3}\left( P_{3}\left(
0,s_{2},s_{3}\right) \right) ds_{2}ds_{3}=-s_{10}I_{2}^{(s)}\left( \mathbf{a}%
_{12},\mathbf{a}_{22},\mathbf{e}_{22}\right) ,  \\
&I_{23}=-s_{20}\int_{0}^{1}\int_{0}^{1}F_{3}\left( P_{3}\left(
s_{1},0,s_{3}\right) \right) ds_{1}ds_{3}=-s_{20}I_{2}^{(s)}\left( \mathbf{a}%
_{13},\mathbf{a}_{23},\mathbf{e}_{23}\right) , \\
&I_{24}=\left( 1+s_{30}\right) \int_{0}^{1}\int_{0}^{1-s_{1}}F_{3}\left(
P_{3}\left( s_{1},s_{2},1\right) \right) ds_{2}ds_{1}=\left( 1+s_{30}\right)
I_{2}^{(t)}\left( \mathbf{a}_{14},\mathbf{a}_{24},\mathbf{e}_{24}\right) , \\
&I_{25}=-s_{30}\int_{0}^{1}\int_{0}^{1-s_{1}}F_{3}\left( P_{3}\left(
s_{1},s_{2},0\right) \right) ds_{2}ds_{1}=-s_{30}I_{2}^{(t)}\left( \mathbf{a}%
_{15},\mathbf{a}_{25},\mathbf{e}_{25}\right) .
\end{aligned}\end{equation}%
We use standardization of 2D domains, which are square (superscript
``$s$'') and triangle (superscript ``$t$''),%
\begin{equation}\begin{aligned}
I_{2}^{(s)}\left( \mathbf{a}_{1j},\mathbf{a}_{2j},\mathbf{e}_{2j}\right) 
&=\int_{D_{2}^{(s)}}...=\int_{0}^{1}\int_{0}^{1}G_{2}\left( R_{2}\left(
s_{1},s_{2};\mathbf{a}_{1j},\mathbf{a}_{2j},\mathbf{e}_{2j}\right) \right)
ds_{2}ds_{1},  \label{b6} \\
I_{2}^{(t)}\left( \mathbf{a}_{1j},\mathbf{a}_{2j},\mathbf{e}_{2j}\right) 
&=\int_{D_{2}^{(t)}}...=\int_{0}^{1}\int_{0}^{1-s_{1}}G_{2}\left(
R_{2}\left( s_{1},s_{2};\mathbf{a}_{1j},\mathbf{a}_{2j},\mathbf{e}%
_{2j}\right) \right) ds_{2}ds_{1}.
\end{aligned}\end{equation}%
Note that $s_{10},s_{20},s_{30}$, and $h_{3}$ can be found from parameters
of $I_{3}$ in (\ref{b5}) (see Appendix), while parameters of the 2D
integrals are related to these parameters as 
\begin{equation}\begin{aligned}
\mathbf{a}_{11} &=\mathbf{a}_{1}-\mathbf{a}_{2},\quad \mathbf{a}_{21}=\mathbf{a}_{3},\quad \mathbf{a}_{12}=\mathbf{a}_{2},\quad \mathbf{a}_{22}=\mathbf{a}_{3},\quad \mathbf{a}_{13}=\mathbf{a}_{1},\quad \mathbf{a}_{23}=\mathbf{a}_{3}, \label{b7} \\
\mathbf{a}_{14} &=\mathbf{a}_{1},\quad \mathbf{a}_{24}=\mathbf{a}_{2},\quad \mathbf{a}_{15}=\mathbf{a}_{1},\quad \mathbf{a}_{25}=\mathbf{a}_{2},\quad  \mathbf{e}_{21}=\mathbf{e}_{22}+\mathbf{a}_{2}, \\
\mathbf{e}_{22} &=\mathbf{e}_{23}=\mathbf{e}_{25}=\mathbf{a}_{1}s_{10}+\mathbf{a}_{2}s_{20}+\mathbf{a}_{3}s_{30},\quad \mathbf{e}_{24}=\mathbf{e}_{22}+\mathbf{a}_{3}.\quad  
\end{aligned}\end{equation}%

\subsubsection{2D domain boundaries and 1D integrals}

Following the scheme of the previous subsections we can reduce integration
over the square and triangle to contour integration consisting of 4 and
3 line integrals, respectively. We have then%
\begin{gather}
I_{2}^{(s)}\left( \mathbf{a}_{1},\mathbf{a}_{2},\mathbf{e}_{2}\right)
=\sum_{j=1}^{4}I_{1j}^{(s)},\quad I_{2}^{(t)}\left( \mathbf{a}_{1},\mathbf{a}%
_{2},\mathbf{e}_{2}\right) =\sum_{j=1}^{3}I_{1j}^{(t)},  \label{b8} \\
I_{11}^{(s)}=\left( 1+s_{10}\right) \int_{0}^{1}F_{2}\left( P_{2}\left(
1,s_{2}\right) \right) ds_{2}=\left( 1+s_{10}\right) I_{1}\left( \mathbf{a}%
_{11},\mathbf{e}_{11}\right) ,  \notag \\
I_{12}^{(s)}=I_{11}^{(t)}=-s_{10}\int_{0}^{1}F_{2}\left( P_{2}\left(
0,s_{2}\right) \right) ds_{2}=-s_{10}I_{1}\left( \mathbf{a}_{12},\mathbf{e}%
_{12}\right) ,  \notag \\
I_{13}^{(s)}=\left( 1+s_{20}\right) \int_{0}^{1}F_{2}\left( P_{2}\left(
s_{1},1\right) \right) ds_{1}=\left( 1+s_{20}\right) I_{1}\left( \mathbf{a}%
_{13},\mathbf{e}_{13}\right) ,  \notag \\
I_{14}^{(s)}=I_{12}^{(t)}=-s_{20}\int_{0}^{1}F_{2}\left( P_{2}\left(
s_{1},0\right) \right) ds_{1}=-s_{20}I_{1}\left( \mathbf{a}_{14},\mathbf{e}%
_{14}\right) ,  \notag \\
I_{13}^{(t)}=\left( 1+s_{10}+s_{20}\right) \int_{0}^{1}F_{2}\left(
P_{2}\left( s_{1},1-s_{1}\right) \right) ds_{1}=\left(
1+s_{10}+s_{20}\right) I_{1}\left( \mathbf{a}_{15},\mathbf{e}_{15}\right) ,
\notag\\
I_{1}\left( \mathbf{a}_{1j},\mathbf{e}_{1j}\right)
=\int_{D_{1}}...=\int_{0}^{1}G_{1}\left( R_{1}\left( s_{1};\mathbf{a}_{1j},%
\mathbf{e}_{1j}\right) \right) ds_{1}.\notag   
\end{gather}%
Also, $s_{10},s_{20},$ and $h_{2}$ can be found from parameters of $I_{2}$
in (\ref{b8}) (see Appendix), while parameters of the 1D integrals are
related to these parameters as%
\begin{eqnarray}
\mathbf{a}_{11} &=&\mathbf{a}_{12}=\mathbf{a}_{2}\mathbf{,\quad a}_{13}=%
\mathbf{a}_{14}=\mathbf{a}_{1},\quad \mathbf{a}_{15}=\mathbf{a}_{1}-\mathbf{a%
}_{2}\mathbf{,}  \label{b10} \\
\mathbf{e}_{12} &=&\mathbf{e}_{14}=\mathbf{a}_{1}s_{10}+\mathbf{a}%
_{2}s_{20},\quad \mathbf{e}_{11}=\mathbf{e}_{12}+\mathbf{a}_{1},\quad 
\mathbf{e}_{13}=\mathbf{e}_{15}=\mathbf{e}_{12}+\mathbf{a}_{2}.\quad   \notag
\end{eqnarray}

\subsubsection{Evaluation of 1D integrals}

The integration process ends with evaluation of the 1D integral (also, see (%
\ref{mp19}))
\begin{equation}
I_{1}\left( \mathbf{a}_{1},\mathbf{e}_{1}\right)  = \left( 1+s_{10}\right)
F_{1}\left( P_{1}\left( 1\right) \right) -s_{10}F_{1}\left( P_{1}\left(
0\right) \right) ,  \label{b11} \quad
P_{1}\left( s_{1}\right)  = \left| s_{1}+s_{10}\right| \left| \mathbf{a}%
_{1}\right| ,  
\end{equation}%
where $F_{1}\left( P_{1}\right) $ is provided by expressions (\ref{k6}) and
parameters $s_{10}$ and $h_{1}$ can be found from $\mathbf{a}_{1}$ and $%
\mathbf{e}_{1}$ as provided in the Appendix.

\section{Gradient of single layer and double layer potential}
We next consider integrals (\ref{ps2}) and (\ref{ps3}). We start
with the case when $S_{x}$ and $S_{y}$ do not have common points, 
$S_{x}\cap S_{y}=\emptyset $. In this case integral $M$ does not require a special treatment, 
since $M=$ $-\mathbf{n}_{x}\cdot \mathbf{L}^{\prime }$. Indeed, due to the symmetry property of $G\left( \mathbf{x,y}\right) $ (\ref{ps6}), 
we have%
\begin{equation}
\mathbf{L}^{\prime }=\int_{S_{y}}\int_{S_{x}}\nabla _{\mathbf{y}}G\left( 
\mathbf{x,y}\right) dS\left( \mathbf{x}\right) dS\left( \mathbf{y}\right)
=-\int_{S_{y}}\int_{S_{x}}\nabla _{\mathbf{x}}G\left( \mathbf{x,y}\right)
dS\left( \mathbf{x}\right) dS\left( \mathbf{y}\right) .  \label{db2}
\end{equation}
Note that the gradient can be decomposed into the surface gradient
and normal derivative operators:%
\begin{equation}
\nabla _{\mathbf{x}}=\nabla _{s\mathbf{x}}+\mathbf{n}_{x}\left( \mathbf{n}%
_{x}\cdot \nabla _{\mathbf{x}}\right) ,\quad \nabla _{\mathbf{y}}=\nabla _{s%
\mathbf{y}}+\mathbf{n}_{y}\left( \mathbf{n}_{y}\cdot \nabla _{\mathbf{y}%
}\right) .  \label{db3}
\end{equation}%
This shows that two forms of the integral in (\ref{db2}) can be written as%
\begin{eqnarray}
\mathbf{L}^{\prime } &=&\mathbf{F}_{y}+\mathbf{n}_{y}L_{y}^{\prime },\quad 
\mathbf{L}^{\prime }=-\mathbf{F}_{x}+\mathbf{n}_{x}L_{x}^{\prime },\quad
\label{db4} \\
L_{y}^{\prime } &=&\mathbf{n}_{y}\cdot \mathbf{L}^{\prime },\quad \mathbf{F}%
_{y}=\int_{S_{x}}\int_{S_{y}}\nabla _{s\mathbf{y}}G\left( \mathbf{x,y}%
\right) dS\left( \mathbf{y}\right) dS\left( \mathbf{x}\right) ,  \notag \\
L_{x}^{\prime } &=&\mathbf{n}_{x}\cdot \mathbf{L}^{\prime },\quad \mathbf{F}%
_{x}=\int_{S_{y}}\int_{S_{x}}\nabla _{s\mathbf{x}}G\left( \mathbf{x,y}%
\right) dS\left( \mathbf{x}\right) dS\left( \mathbf{y}\right) .  \notag
\end{eqnarray}%
Taking the scalar product of the first equation with $\mathbf{n}_{x}$ and
the second equation with $\mathbf{n}_{y}$, we obtain simultaneous
equations with respect to $L_{x}^{\prime }$ and $L_{y}^{\prime },$ 
\begin{equation}
L_{x}^{\prime } = \mathbf{n}_{x}\cdot \mathbf{F}_{y}+\left( \mathbf{n}%
_{x}\cdot \mathbf{n}_{y}\right) L_{y}^{\prime },  \quad 
L_{y}^{\prime } = -\mathbf{n}_{y}\cdot \mathbf{F}_{x}+\left( \mathbf{n}%
_{x}\cdot \mathbf{n}_{y}\right) L_{x}^{\prime }.  \label{db5} \\
\end{equation}%
This system is non-degenerate unless the planes $P_{x}$ and $P_{y}$
containing triangles $S_{x}$ and $S_{y}$ are parallel ($\left| \mathbf{n}%
_{x}\cdot \mathbf{n}_{y}\right| =1$). Below we consider both, the
non-degenerate and degenerate cases. Before that we note that integrals $%
\mathbf{F}_{x}$ and $\mathbf{F}_{y}$ can be reduced to 3D integrals using
the Gauss' theorem, 
\begin{equation}\begin{aligned}\label{db7} 
\mathbf{F}_{x} =&\int_{S_{y}}\int_{C_{x}}\mathbf{n}_{cx}G\left( \mathbf{x,y}%
\right) dC\left( \mathbf{x}\right) dS\left( \mathbf{y}\right) =\sum_{i=1}^{3}%
\mathbf{n}_{cxi}F_{xi},\\ 
F_{xi} =& \int_{S_{y}}\int_{C_{xi}}G\left( \mathbf{x,y}\right) dC\left( \mathbf{x}\right) dS\left( \mathbf{y}\right) ,\quad 
\mathbf{n}_{cxi} = \frac{\left( \mathbf{l}_{xi}\times \mathbf{n}_{x}\right)}{l_{xi}},\quad 
\mathbf{n}_{x} = -\frac{ \mathbf{l}_{x1}\times \mathbf{l}_{x3} }{2A_{x}} ,\\
\mathbf{F}_{y} =&\int_{S_{x}}\int_{C_{y}}\mathbf{n}_{cy}G\left( \mathbf{x,y}%
\right) dC\left( \mathbf{y}\right) dS\left( \mathbf{x}\right) =\sum_{i=1}^{3} \mathbf{n}_{cyi}F_{yi},\\
F_{yi}=&\int_{S_{x}}\int_{C_{yi}}G\left( \mathbf{%
x,y}\right) dC\left( \mathbf{y}\right) dS\left( \mathbf{x}\right) , \quad
  \mathbf{n}_{cyi} = \frac{ \mathbf{l}_{yi}\times \mathbf{n}_{y} }{l_{yi}},\quad 
  \mathbf{n}_{y} = -\frac{ \mathbf{l}_{y1}\times \mathbf{l}_{y3} }{2A_{y}} .  
\end{aligned}\end{equation}%
Here $\mathbf{n}_{cx}$ and $\mathbf{n}_{cy}$ are the outer normals to
contours $C_{x}$ and $C_{y}$ bounding triangles $S_{x}$ and $S_{y}$,
respectively. These contours consist of segments $C_{xi}$ and $C_{yi}$ with
normals $\mathbf{n}_{cxi}$ and $\mathbf{n}_{cyi}$, $i=1,2,3.$

So, we have six scalar integrals (\ref{db7}), which are nothing, but
integrals over prismatic domains (\ref{b1}) where%
\begin{equation}
F_{yi} = 2l_{yi}A_{x}\widetilde{I}_{3}\left( \mathbf{a}_{1i},\mathbf{a}%
_{2i},\mathbf{a}_{3i},\mathbf{e}_{3i}\right) ,\quad F_{xi} = 2l_{xi}A_{y}\widetilde{I}_{3}\left( \mathbf{a}_{1j},\mathbf{a%
}_{2j},\mathbf{a}_{3j},\mathbf{e}_{3j}\right) ,
\label{db7.1}
\end{equation}%
for $i=1,2,3$, $j=i+3$, and $\widetilde{I}_{3}$ is defined as 
\begin{equation}\begin{aligned}
&\widetilde{I}_{3}\left( \mathbf{a}_{1j},\mathbf{a}_{2j},\mathbf{a}_{3j},\mathbf{e}_{3j}\right)
=\int_{D_{3}}...= \int_{0}^{1}\int_{0}^{1-s_{1}}\int_{0}^{1}\widetilde{G}_{3}\left( R_{3}\left( s_{1},s_{2},s_{3};\mathbf{a}_{1j},\mathbf{a}_{2j},\mathbf{a}_{3j},\mathbf{e}_{3j}\right) \right) ds_{3}ds_{2}ds_{1},
\label{db7.2} \\
&\widetilde{G}_{3}\left( R_{3}\right) = \left(
R_{3}^{2}+h_{4}^{2}\right) ^{-1/2}=G_{4}\left( R_{4}\right) =R_{4}^{-1}, 
\end{aligned}\end{equation}%
and parameters $\mathbf{a}_{1j},\mathbf{a}_{2j},\mathbf{a}_{3j},\mathbf{e}%
_{3j}$, $j=1,...,6$ are provided by (\ref{b4}), where the values for
parameters $\mathbf{a}_{1},\mathbf{a}_{2},\mathbf{a}_{3},\mathbf{a}_{4}$ of
the 4D integrals are provided by (\ref{b01}) and coordinates of projections $%
s_{10},...,s_{40}$ and $h_{4}$ can be found from these parameters and $%
\mathbf{e}_{4}$ in the same way as for the single layer potential. So, the
integrals $I_{3}$ and $\widetilde{I}_{3}$ depend on the same parameter set,
but have different integrands.

\subsection{Non-degenerate case}

Solving the system for the non-degenerate case we get 
\begin{equation}\begin{aligned}
-M = L_{x}^{\prime }=\frac{\mathbf{n}_{x}\cdot \mathbf{F}_{y}-\left( 
\mathbf{n}_{x}\cdot \mathbf{n}_{y}\right) \mathbf{n}_{y}\cdot \mathbf{F}_{x}%
}{1-\left( \mathbf{n}_{x}\cdot \mathbf{n}_{y}\right) ^{2}}, \quad \label{db6} 
\mathbf{L}^{\prime } = -\mathbf{F}_{x}+\mathbf{n}_{x}L_{x}^{\prime }=-%
\mathbf{F}_{x}-\mathbf{n}_{x}M. 
\end{aligned}\end{equation}%
This shows that it is sufficient to compute integrals $\mathbf{F}_{x}$ and $%
\mathbf{F}_{y}$ to obtain the values of integrals 
(\ref{ps2}) - (\ref{ps3}) 
in the
non-degenerate case. 
The non-degenerate case corresponds to $\left| \mathbf{n}_{x}\cdot \mathbf{n}%
_{y}\right| \neq 1$, so $P_{x}$ and $P_{y}$ are not parallel and $h_{4}=0$.
Note that in this case we have (see (\ref{k3}) and (\ref{b3}))%
\begin{equation}
\widetilde{G}_{3}\left( R_{3}\right) =3G_{3}\left( R_{3}\right) ,\quad 
\widetilde{I}_{3}\left( \mathbf{a}_{1j},\mathbf{a}_{2j},\mathbf{a}_{3j},%
\mathbf{e}_{3j}\right) =3I_{3}\left( \mathbf{a}_{1j},\mathbf{a}_{2j},\mathbf{%
a}_{3j},\mathbf{e}_{3j}\right) .  \label{db7.3}
\end{equation}%
Hence, if one needs to compute integrals (\ref{ps1}) - (\ref{ps3}) for the
most common situation $h_{4}=0$, integrals (\ref{ps2}) - (\ref{ps3}) come
for almost no additional cost.

\subsection{Degenerate case}

Note that after $M$ is computed $\mathbf{L}^{\prime }$ can be found from the
last equation (\ref{db6}. There are two problems that need to be resolved for the degenerate case. 
First, how to compute $M$ as the method presented by (\ref{db6}) does not work in this case. 
Second, how to evaluate integrals (\ref{db7.2}) at $h_{4}\neq 0$.
\subsubsection{Evaluation of $M$}

In case when planes $P_{x}$ and $P_{y}$ are parallel we can consider the function%
\begin{equation}
G^{\prime }\left( \mathbf{x,y}\right) =\mathbf{n}_{x}\cdot \nabla _{\mathbf{x%
}}G\left( \mathbf{x},\mathbf{y}\right) =\frac{\mathbf{n}_{x}\cdot \left( 
\mathbf{y}-\mathbf{x}\right) }{\left| \mathbf{x}-\mathbf{y}\right| ^{3}}=%
\frac{\delta }{\left| \mathbf{x}-\mathbf{y}\right| ^{3}},\quad h_{4}=\left|
\delta \right|  \label{s8}
\end{equation}%
where constant $\delta $ is a signed distance between the planes of $P_{x}$
and $P_{y}$. So, 
\begin{gather}
M =\int_{S_{y}}\int_{S_{x}}G^{\prime }\left( \mathbf{x,y}\right) dS\left( 
\mathbf{x}\right) dS\left( \mathbf{y}\right) =4A_{x}A_{y}\delta U^{\prime },
\label{s9} \\ \notag
U^{\prime } =\int_{0}^{1}\!\!\int_{0}^{1-u}\!\!\!\int_{0}^{1}\!\!\int_{0}^{1-s}\!\!\!\!\! g^{\prime
}\left( s,t,u,v\right) dtdsdvdu,\quad g^{\prime }=\frac{G^{\prime
}\left( \mathbf{X}\left( s,t\right) \mathbf{,Y}\left( u,v\right) \right)}{\delta }
= \frac{1}{\left| \mathbf{X}-\mathbf{Y}\right| ^{3}}. 
\end{gather}%
Here and below we use the prime for the integrals related to $M$. It is also
clear that at $\delta =0$ we have $G^{\prime }=0$ ($S_{x}$ and $S_{y}$ in
the same plane). In fact, this corresponds to the principal value of
integral (\ref{ps2}) when $\mathbf{y}\in S_{x}$. So, we should set $M=0$ for 
$\delta =0$.

Now computations of $U^{\prime }$ are similar to computations of $U$ given
by (\ref{s6}) with the only difference that instead of (\ref{k1}) we should
use PBF%
\begin{equation}
G_{4}^{\prime }\left( R_{4}\right) =R_{4}^{-3}  \label{db8}
\end{equation}%
and consider only cases \#6 and \#7 listed in (\ref{b0.1}). 
The results of integration for $d=4,3,$ and $2$ are: 
\begin{equation}\begin{aligned}
\#6,7 &:\quad F_{4}^{\prime }\left( P_{4};h_{4}\right) =\frac{1}{P_{4}^{4}%
}\left( \frac{P_{4}^{2}+2h_{4}^{2}}{R_{4}}-2h_{4}\right) ,\quad R_{4}=\sqrt{%
P_{4}^{2}+h_{4}^{2}},  \label{db9} \\
\#6,7 &:\quad F_{3}^{\prime }(P_{3};0,h_{4})=\frac{1}{P_{3}^{3}}\left[ -%
\frac{2P_{3}}{R_{4}+h_{4}}+\ln \frac{P_{3}+R_{4}}{h_{4}}\right] ,\quad R_{4}=%
\sqrt{P_{3}^{2}+h_{4}^{2}},  \\
\#6 &:\quad\! F_{2}^{\prime }\!\left( P_{2};0,0,h_{4}\right) =\frac{1}{P_{2}^{2}%
}\!\left[ \frac{1}{R_{4}+h_{4}}-\frac{1}{P_{2}}\ln \frac{P_{2}+R_{4}}{h_{4}}+%
\frac{1}{2h_{4}}\right] ,\quad\! R_{4}=\sqrt{P_{2}^{2}+h_{4}^{2}}, \\
\#7 &:\quad\! F_{2}^{\prime }\!\left( P_{2};h_{2},0,h_{4}\right)= \\&\!\! \frac{1}{%
P_{2}^{2}}\!\left( \frac{1}{R_{4}+h_{4}}-\frac{1}{h+h_{4}}-\frac{1}{R_{2}}\ln 
\frac{R_{2}+R_{4}}{h_{4}}+\frac{1}{h_{2}}\ln \frac{h_{2}+h}{h_{4}}\right)
,\quad R_{4}=\sqrt{P_{2}^{2}+h^{2}}, 
\end{aligned}\end{equation}%
with $h=\sqrt{h_{2}^{2}+h_{4}^{2}}$. For $d=1$ we have 
\begin{equation}\begin{aligned}
\#6 :&\quad\!\!\! F_{1}^{\prime }\!\left( P_{1};h_{1},0,0,h_{4}\right) =\frac{1}{%
h_{1}^{2}}\!\!\left[ \Phi_{1}\! \left( P_{1}\right) +\frac{h_{1}^{2}-h_{4}^{2}}{%
2h_{1}h_{4}}\Phi _{2}\!\left( P_{1};h_{1}\right) -\Phi _{3}\!\left( P_{1}\right)
+\frac{1}{2\left( R_{4}+h_{4}\right) }\right] ,  \label{db9.1} \\
\#7 :&\quad\!\!\! F_{1}^{\prime }\!\left( P_{1};0,h_{2},0,h_{4}\right) =\frac{-1}{%
h_{2}^{2}}\left[ \Phi _{1}\left( P_{1}\right) - \frac{h_{4}}{h_{2}}\Phi
_{2}\left( P_{1};h_{2}\right) - \frac{h_{2}^{2}}{h_{4}^{2}}\Phi _{4}\left(
P_{1};h_{4}\right) + \frac{1}{P_{1}^{2}}\left( \frac{R_{2}^{2}}{R_{4}+h_{4}}-%
\frac{h_{2}^{2}}{h+h_{4}}\right) \right] , 
\end{aligned}\end{equation}%
where $R_{2},R_{4},\Phi _{1},\Phi _{2},\Phi _{3},$ and $\Phi _{4}$ are
provided by (\ref{k6}).

\subsubsection{Evaluation of $\mathbf{F}_{x}$}

Again, we consider only cases \#6 and \#7 from (\ref{b0.1}) for the PBFs
generated by $\widetilde{G}_{3}\left( R_{3}\right) $, which we mark with
tilde. With $R_{4} = \sqrt{P_{3}^{2}+h_{4}^{2}}$, we have by definition
\begin{equation}\begin{aligned}
\#6,7: \quad\!\!\! \widetilde{F}_{3}\left( P_{3};0,h_{4}\right) =& \frac{1}{P_{3}^{3}%
}\int_{0}^{P_{3}}\!\!\!\! P^{2}\widetilde{G}_{3}\left( \sqrt{P^{2}+h_{4}^{2}}\right)
dP = \frac{R_{4}}{2P_{3}^2}-\frac{h_{4}^{2}}{2P_3^3}\ln \frac{P_{3}+R_{4}}{%
h_{4}}.  \label{db9.2}
\end{aligned}\end{equation}%
Comparing this with (\ref{k4}) we can see that this function
coincides with $F_{3}\left( P_{3};h_{3},0\right) $ multiplied by 3, but with
parameter $h_{4}$ instead of $h_{3}$, which means
\begin{equation}
\widetilde{F}_{3}\left( P_{3};0,h_{4}\right) =3F_{3}\left(
P_{3};h_{4},0\right) .  \label{db9.3}
\end{equation}%
Consequently, at $d=2$ and $d=1$ cases \#6 and \#7 for $\widetilde{F}_{d}$
can be mapped to cases \#4 and \#5 for $F_{d}$. So, 
\begin{equation}\begin{aligned}
&\#6\quad\! \widetilde{F}_{2}\left( P_{2};0,0,h_{4}\right) = 3F_{2}\left(
P_{2};0,h_{4},0\right) ,\quad\! \widetilde{F}_{1}\left(
P_{1};h_{1},0,0,h_{4}\right) =3F_{1}\left( P_{1};h_{1},0,h_{4},0\right) ,
\label{db9.4} \\
&\#7\quad\! \widetilde{F}_{2}\!\left( P_{2};h_{2},0,h_{4}\right) = 3F_{2}\!\left(
P_{2};h_{2},h_{4},0\right) ,\quad \!
\widetilde{F}_{1}\!\left(
P_{1};0,h_{2},0,h_{4}\right) =3F_{1}\!\left( P_{1};0,h_{2},h_{4},0\right) . 
\end{aligned}\end{equation}

\section{Normal derivative \textbf{of double layer}}

Consider now integrals (\ref{ps4}) 
in a non-singular case, $%
S_{x}\cap S_{y}=\emptyset $. Let us introduce vector function 
\begin{equation}
\mathbf{u}\left( \mathbf{y}\right) =\mathbf{n}_{x}\int_{S_{x}}G\left( 
\mathbf{x,y}\right) dS\left( \mathbf{x}\right) .  \label{gdb2}
\end{equation}%
We have then 
\begin{equation}
\nabla _{\mathbf{y}}\cdot \mathbf{u}\left( \mathbf{y}\right) =\mathbf{n}%
_{x}\cdot \nabla _{\mathbf{y}}\int_{S_{x}}G\left( \mathbf{x,y}\right)
dS\left( \mathbf{x}\right) =-\int_{S_{x}}\mathbf{n}_{x}\cdot \nabla _{%
\mathbf{x}}G\left( \mathbf{x,y}\right) dS\left( \mathbf{x}\right) .
\label{gdb3}
\end{equation}%
Furthermore, as follows from definition (\ref{gdb2}) $\mathbf{u}\left( 
\mathbf{y}\right) $ is a harmonic function ($\nabla _{\mathbf{y}}^{2}\mathbf{%
u}=\mathbf{0}$). So, we have 
\begin{equation}
\nabla _{\mathbf{y}}\left[ \nabla _{\mathbf{y}}\cdot \mathbf{u}\left( 
\mathbf{y}\right) \right] =\nabla _{\mathbf{y}}^{2}\mathbf{u}+\nabla _{%
\mathbf{y}}\times \nabla _{\mathbf{y}}\times \mathbf{u}\left( \mathbf{y}%
\right) =\nabla _{\mathbf{y}}\times \nabla _{\mathbf{y}}\times \mathbf{u}%
\left( \mathbf{y}\right) .  \label{gdb4}
\end{equation}%
Relations (\ref{gdb3}) and (\ref{gdb4}) show then that integral (\ref{ps4})
can be represented in the form%
\begin{equation}\begin{aligned}
M^{\prime } =& -\int_{S_{y}}\mathbf{n}_{y}\cdot \nabla _{\mathbf{y}}\times
\int_{S_{x}}\mathbf{n}_{x}\times \nabla _{\mathbf{x}}G\left( \mathbf{x,y}%
\right) dS\left( \mathbf{x}\right) dS\left( \mathbf{y}\right) ,\quad
M^{\prime }=\mathbf{n}_{y}\cdot \mathbf{M}^{\prime },\\  \label{gdb5}
\mathbf{M}^{\prime } =& \int_{S_{y}} \nabla_{\mathbf{y}} \int_{S_{x}} \mathbf{n}_{x} \cdot \nabla_{\mathbf{x}} G\left( \mathbf{x,y}%
\right) dS\left( \mathbf{x}\right) dS\left( \mathbf{y}\right).\\
\end{aligned}\end{equation}%
Here, we put $\nabla _{\mathbf{y}}\times $ under integral over $S_{x}$ and
used identity 
\begin{equation}
\nabla _{\mathbf{y}}\times \left[ \mathbf{n}_{x}G\left( \mathbf{x,y}\right) %
\right] =-\mathbf{n}_{x}\times \nabla _{\mathbf{y}}G\left( \mathbf{x,y}%
\right) =\mathbf{n}_{x}\times \nabla _{\mathbf{x}}G\left( \mathbf{x,y}%
\right) .  \label{gdb6}
\end{equation}%
The surface integral over $S_{x}$ in (\ref{gdb5}) can be reduced to the
contour integral,%
\begin{equation}
\int_{S_{x}}\mathbf{n}_{x}\times \nabla _{\mathbf{x}}G\left( \mathbf{x,y}%
\right) dS\left( \mathbf{x}\right) =\int_{C_{x}}G\left( \mathbf{x,y}\right) d%
\mathbf{x.}  \label{gdb7}
\end{equation}%
Integral $M^{\prime }$ then can be calculated using the Stokes' theorem,%
\begin{eqnarray}
M^{\prime } &=&-\int_{S_{y}}\mathbf{n}_{y}\cdot \nabla _{\mathbf{y}}\times
\int_{C_{x}}G\left( \mathbf{x,y}\right) d\mathbf{x}dS(\mathbf{y}%
)=-\int_{C_{y}}\int_{C_{x}}G\left( \mathbf{x,y}\right) d\mathbf{x}\cdot d%
\mathbf{y}  \label{gdb8} \\
&=&-\sum_{i=1}^{3}\sum_{j=1}^{3}\frac{\mathbf{l}_{xi}\cdot \mathbf{l}_{yj}}{%
l_{xi}l_{yj}}H_{ij},\text{\quad }H_{ij}=\int_{C_{xi}}\int_{C_{yj}}G\left( 
\mathbf{x,y}\right) dC\left( \mathbf{y}\right) dC\left( \mathbf{x}\right) . 
\notag
\end{eqnarray}%
Here we decomposed the integration into integration along the edges $C_{xj}$
and $C_{yi}$ of triangles $S_{x}$ and $S_{y}$ which have lengths $l_{xj}$
and $l_{yi}$ are directed along unit vectors $\mathbf{l}_{xj}$/ $l_{xj}$ and 
$\mathbf{l}_{yi}/l_{yi}$ , respectively.

Obviously, integrals $H_{ij}$ can be expressed as integrals over a square, 
\begin{equation}\begin{aligned}
&H_{ij} = l_{xi}l_{yj}\int_{0}^{1}\int_{0}^{1}\frac{ds_{2}ds_{1}}{\left| 
\mathbf{l}_{xi}s_{1}-\mathbf{l}_{yj}s_{2}+\mathbf{x}_{i}-\mathbf{y}%
_{j}\right| }=l_{xi}l_{yj}\widehat{I}_{2}^{(s)}\left( \mathbf{l}_{xi},-%
\mathbf{l}_{yj},\mathbf{x}_{i}-\mathbf{y}_{j}\right) ,
\label{gdb11} \\
&\widehat{I}_{2}^{(s)}\left( \mathbf{a}_{1},\mathbf{a}_{2},\mathbf{e}%
_{2}\right) = \int_{0}^{1}\int_{0}^{1}\widehat{G}_{2}\left( \widehat{R}%
_{2}\left( s_{1},s_{2};\mathbf{a}_{1},\mathbf{a}_{2},\mathbf{e}_{2}\right)
\right) ds_{2}ds_{1},\\
& \widehat{G}_{2}\left( \widehat{R}_{2}\right) =%
\widehat{R}_{2}^{-1},\quad \widehat{R}_{2}=\left| \mathbf{a}_{1}s_{1}+%
\mathbf{a}_{2}s_{2}+\mathbf{e}_{2}\right| , \quad i,j = 1,2,3.
\end{aligned}\end{equation} 
Here and below the hat is used to mark integrals with integrand produced by $%
\widehat{G}_{2}\left( \widehat{R}_{2}\right) $. Using (\ref{mp14}) we find 
\begin{equation}\begin{aligned}
\widehat{F}_{2}\left( P_{2};h_{2},0,0\right) =&\frac{1}{P_{2}^{2}}\int_{0}^{%
\widehat{P}_{2}}\frac{PdP}{\sqrt{P^{2}+h_{2}^{2}}} =\frac{1}{P_{2}^{2}}%
\left( R_{2}-h_{2}\right) =6F_{2}\left( P_{2};h_{2},0,0\right), \\
\widehat{F}_{1}\left( P_{1};h_{1},h_{2},0,0\right) =& 6F_{1}\left(
P_{2};h_{1},h_{2},0,0\right),  \quad R_{2}=%
\sqrt{P_{2}^{2}+h_{2}^{2}},  \label{gdb12} 
\end{aligned}\end{equation}%
where $F_{2}$ and $F_{1}$ should be taken from (\ref{k5}) and (\ref{k6}) for
cases \#1, \#2, and \#3. From the relation between (\ref{gdb12}) and (\ref%
{k5}) it is clear that at $h_{4}=h_{3}=0$ the integrand in 2D integrals is simply 
$G(R_{2})/6.$ The only remark here is that $h_{1}$ and $h_{2}$ in (\ref%
{gdb12}) should be obtained based on decomposition of $\mathbf{e}_{2}$
appearing in (\ref{gdb11}) and related to the difference of the coordinates
of vertices of triangles $S_{x}$ and $S_{y}$.

Note that a special case $h_{1}=h_{2}$ $=0$ (case \#8)
appears for integrals (\ref{gdb11}). This happens when $\mathbf{l}_{xi},%
\mathbf{l}_{yj}$, and $\mathbf{e}_{2}=\mathbf{x}_{i}-\mathbf{y}_{j}$ are all
collinear vectors. Indeed, since $\mathbf{e}_{2}\in span\left( \mathbf{l}%
_{xi},\mathbf{l}_{yj}\right) $ we have $h_{2}=0$. Furthermore, since $%
\mathbf{l}_{xi}\times \mathbf{l}_{yj}=0$, \ we have $\mathbf{e}_{2}=\mathbf{e%
}_{2\bot }=\mu \mathbf{l}_{xi}$. According (\ref{b8})-(\ref{b10}) we have then 
\begin{equation}
\widehat{I}_{2}^{(s)}\!\left( \mathbf{l}_{xi},-\mathbf{l}_{yj},\mu \mathbf{l}%
_{xi}\right) =\left( 1+\mu \right) \widehat{I}_{1}\!\left( -\mathbf{l}%
_{yj},\left( 1+\mu \right) \mathbf{l}_{xi}\right) -\mu \widehat{I}_{1}\!\left(
-\mathbf{l}_{yj},\mu \mathbf{l}_{xi}\right) +\widehat{I}_{1}\left( \mathbf{l}%
_{xi},\mu \mathbf{l}_{xi}-\mathbf{l}_{yj}\right) .  \label{gdb14}
\end{equation}%
Here again $\mathbf{e}_{1}$ is collinear to $\mathbf{a}_{1}$ for each 1D
integral. This means $h_{1}=0$. \ Integrals for case \#8 can have removable
singularities or diverge. Indeed, in this case form (\ref{gdb11}) shows that
if the denominator turns to zero at some internal point $0<s_{1}<1$ and $%
0<s_{2}<1$ the integral diverges, while if the denominator is not zero at $%
0<s_{1}<1$ and $0<s_{2}<1$ the integrals exist.

\section{Treatment of singularities}

All integrals (\ref{ps1})-(\ref{ps4}) are singular when $S_{x}\cap S_{y}\neq
\emptyset $ with various degree of singularity and, therefore, some comments
should be made concerning the analytical computation. 
The singularity of (\ref{ps1}) is weak and integration does not have
any problems. 
The singularity of (\ref{ps2}) is strong, but is easily treatable as
for any $\mathbf{y}\in S_{x}$ the principal value of the inner integral is
zero, while this integral is finite for any $\mathbf{y}\in S_{x}^{\pm }$.
So, the integration over $S_{y}$ can be performed in a regular way by
ignoring contributions of points $\mathbf{y}\in S_{x}$. 
The same situation holds with (\ref{ps3}). Despite the inner integral
exists at any $\mathbf{y}$ the normal to $S_{x}$ component of the gradient
is discontinuous at $\mathbf{y}\in S_{x}$, while the components tangential
to $S_{x}$ are continuous. That means that (\ref{ps3}) can be computed
according to the last equation in (\ref{db6}), where the tangential
component $\mathbf{F}_{x}$ can be computed as the contour integral without
any problems. 
An integrable singularity takes place when computing (\ref{ps4}) for
the case when $S_{x}$ and $S_{y}$ have a common vertex. This corresponds to $%
\mathbf{x}_{i}=\mathbf{y}_{j}$, while the other vertices of $S_{x}$ and $%
S_{y}$ are different. This is discussed in a subsection below. 
A non-integrable singularity (``hypersingularity'') manifests itself
when computing (\ref{ps4}) for the case when $S_{x}$ and $S_{y}$ have a
common edge. This corresponds to $\mathbf{l}_{xi}=\mathbf{l}_{yj}$, $\mathbf{%
x}_{i}=\mathbf{y}_{j}$, or $\mathbf{l}_{xi}=-\mathbf{l}_{yj}$, $\mathbf{x}%
_{i}=\mathbf{y}_{j}+\mathbf{l}_{yj}$ in (\ref{gdb11}), at which $H_{ij}$
diverges. Two issues should be considered here. First, the validity of
construction leading to (\ref{gdb8}) for the case $S_{x}\cap S_{y}\neq
\emptyset $, and, second, how this singularity should be handled in practical
computations. We dedicate a subsection below for this discussion.

\subsection{Hypersingular integral}

There exist a number of papers related to the Galerkin hypersingular integrals %
\cite{Salvadori2001CMAME}, \cite{Gray2001EABE}-\cite{Polimeridis2011IEEE}. 
We consider an approach
close to that developed in \cite{Gray2001EABE}, where the singularity is
studied based on the results for a non-singular case $S_{x}\cap
S_{y}=\emptyset $, but the vertex of $S_{y}$, its edge or the entire
triangle $S_{y}$ are located at a small distance from the respective
counterparts of $S_{x}$ and that distance tends to zero. The fact that $%
S_{x}\cap S_{y}=\emptyset $ allows us to use the formalism (\ref{gdb8}) and
focus on investigation of double integrals $H_{ij}$ (\ref{gdb11}). Three
basic cases can then be considered for the case when $S_{y}$ approaches $%
S_{x}$. First, that in the limit $S_{y}$ and $S_{x}$ have a common vertex
(one-touch case), second, that in the limit $S_{y}$ and $S_{x}$ have a
common edge (two-touch case), and, third, that in the limit $S_{y}$ and $%
S_{x}$ coincide (three-touch, or self-action case).

\subsubsection{One-touch case}

Without any loss of generality we can index the vertices in the way that
the common vertex in the limit is $\mathbf{x}_{1}=\mathbf{y}_{1}$. To make $%
S_{x}\cap S_{y}=\emptyset $ we place $\mathbf{y}_{1}$at position%
\begin{equation}
\mathbf{y}_{1}=\mathbf{x}_{1}-\epsilon \mathbf{n}_{11}\notin S_{x},\quad 
\mathbf{n}_{11}\cdot \mathbf{l}_{x1}=0,\quad \mathbf{n}_{11}\cdot \mathbf{l}%
_{y1}=0,\quad \epsilon >0,  \label{si1}
\end{equation}%
We have then for integral (\ref{gdb11})%
\begin{equation}
\widehat{I}_{2}^{(s)}\left( \mathbf{l}_{x1},-\mathbf{l}_{y1},\mathbf{x}_{1}-%
\mathbf{y}_{1}\right) = \widehat{I}_{2}^{(s)}\left( \mathbf{l}_{x1},-%
\mathbf{l}_{y1},\epsilon \mathbf{n}_{11}\right) ,  \label{si2} \quad
s_{10} = 0,\quad s_{20}=0,\quad h_{2}=\epsilon . 
\end{equation}%
Using (\ref{b8}) and (\ref{b10}) we can reduce this to a sum of 1D integrals,%
\begin{equation}
\widehat{I}_{2}^{(s)}\left( \mathbf{l}_{x1},-\mathbf{l}_{y1},\epsilon 
\mathbf{n}_{11}\right) =\widehat{I}_{1}\left( -\mathbf{l}_{y1},\mathbf{l}%
_{x1}\right) +\widehat{I}_{1}\left( \mathbf{l}_{x1},-\mathbf{l}_{y1}\right) .
\label{si3}
\end{equation}%
If $\mathbf{l}_{x1}$ and $\mathbf{l}_{y1}$ are not collinear, for
both integrals $h_{1}\neq 0$. The PBF (\ref{gdb12}) then has a finite limit,
(\ref{k6}) case \#3:%
\begin{equation}
\lim_{h_{2}\rightarrow 0}\widehat{F}_{1}\left( P_{1};h_{1},h_{2},0,0\right) =%
\widehat{F}_{1}\left( P_{1};h_{1},0,0,0\right) =6F_{1}\left(
P_{1};h_{1},0,0,0\right) .  \label{si4}
\end{equation}%
Detailed steps to evaluate integrals in (\ref{si3}) are illustrated in the appendix, and shows that the one-touch case is integrable in a regular way.

It may happen that $\mathbf{l}_{x1}$ and $\mathbf{l}_{y1}$ are collinear. In
this case $h_{1}=0$ and in the limit $h_{2}\rightarrow 0$ a new
limiting case \#8 corresponding to $h_{1}=h_{2}=h_{3}=h_{4}=0$ appears only for hypersingular integrals, 
\begin{equation}
\widehat{F}_{1}\left( P_{1};0,h_{2},0,0\right) =6F_{1}\left(
P_{1};0,h_{2},0,0\right) =\frac{1}{P_{1}}\ln \frac{P_{1}+R_{2}}{h_{2}}-\frac{%
1}{R_{2}+h_{2}}.  \label{si5}
\end{equation}
However, for the one-touch case the finite limit exists. Indeed, for
collinear $\mathbf{l}_{x1}$ and $\mathbf{l}_{y1}$ we have $\mathbf{l}%
_{y1}=\lambda \mathbf{l}_{x1}$, where $\lambda <0$ for the one-touch case
(otherwise $S_{x}$ and $S_{y}$ should have more than one common points). We
have then at $h_{2}\rightarrow 0$ (see (\ref{b11})). 
\begin{eqnarray}
\widehat{I}_{1}\left( -\mathbf{l}_{y1},\mathbf{l}_{x1}\right)  &=&\frac{1}{%
l_{y1}}\lim_{h_{2}\rightarrow 0}\left( \ln \frac{2\left| 1-\frac{1}{\lambda }%
\right| l_{y1}}{h_{2}}-\ln \frac{2\left| \frac{1}{\lambda }\right| l_{y1}}{%
h_{2}}\right) =\frac{1}{l_{y1}}\ln \left( 1-\lambda \right) ,  \label{si5.1}
\\
\widehat{I}_{1}\left( \mathbf{l}_{x1},-\mathbf{l}_{y1}\right)  &=&\frac{1}{%
l_{x1}}\lim_{h_{2}\rightarrow 0}\left( \ln \frac{2\left| 1-\lambda \right|
l_{xi}}{h_{2}}-\ln \frac{2\left| \lambda \right| l_{xi}}{h_{2}}\right) =%
\frac{1}{l_{x1}}\ln \frac{1-\lambda }{\left| \lambda \right| }.  \notag
\end{eqnarray}%
These forms show that a practical way to handle this singularity is to
define the PBF for case \#8 as%
\begin{equation}
\#8:\quad F_{1}\left( P_{1};0,0,0,0\right) =\frac{1}{6P_{1}}\ln P_{1},\quad 
\widehat{F}_{1}\left( P_{1};0,0,0,0\right) =6F_{1}\left(
P_{1};0,0,0,0\right) .\quad   \label{si5.2}
\end{equation}%
Note that this definition is valid only in case $\mathbf{l}_{y1}=\lambda 
\mathbf{l}_{x1}$, $\lambda <0.$

\subsubsection{Two-touch case}

In this case we index the vertices the same way as in the one-touch case,
while having the common edge $\mathbf{l}_{x1}=\mathbf{l}_{y1}$. The $%
\epsilon $ shift of $\mathbf{y}_{1}$ can be performed in the direction
normal to the plane of triangle $S_{x}$. In this case we have the same
equations as (\ref{si1})-(\ref{si5}). The difference is that now in (\ref%
{si3}) the arguments of $\widehat{I}_{1}$ are collinear, $\mathbf{l}%
_{y1}=\lambda \mathbf{l}_{x1},$ $\lambda =1>0$. 
From (\ref{gdb11}), (\ref{si3}) and (\ref{b11}) we have at $\epsilon
\rightarrow 0$%
\begin{equation}
H_{11}=2l_{x1}^{2}\left( \frac{1}{l_{x1}}\ln \frac{l_{x1}+\left(
l_{x1}^{2}+\epsilon ^{2}\right) ^{1/2}}{\epsilon }-\frac{1}{\left(
l_{x1}^{2}+\epsilon ^{2}\right) ^{1/2}+\epsilon }\right) \sim 2l_{x1}\left(
\ln \frac{2l_{x1}}{\epsilon }-1\right).  \label{si6}
\end{equation}%
This illustrates the character of the singularity (logarithmic) of $H_{ij}$
corresponding to the common edge $ij$.

Equation (\ref{gdb8}) shows that in computation of $M^{\prime }$ a factor 
$\mathbf{l}_{xi}\cdot \mathbf{l}_{yj}$ is present. This means
that if we have two different triangles $S_{y1}$ and $S_{y2}$ with a common
edge, $S_{x}=S_{y1}$ and we consider a sum of integrals with the inner
integral over $S_{x}$ and the outer integrals over $S_{y1}$ and $S_{y2}$
then the contributions over the common edge cancel out since $\mathbf{l}%
_{xi}\cdot \mathbf{l}_{yj}=1$ for $S_{y1}$ and $\mathbf{l}_{xi}\cdot \mathbf{%
l}_{yj}=-1$ for $S_{y2}$ (as the contours should be oriented consistently to
provide the same orientation of the normals). 
So, a practical rule here is just zero the contribution over the common edge
for all hypersingular integrals appearing in the Galerkin BEM, 
\begin{equation}
H_{ij}=0,\quad C_{xi}=C_{yj}.  \label{si7}
\end{equation}

\begin{remark}
Formally, the Galerkin BEM with constant panel basis function does not make
sense for hypersingular integral equations. Indeed, in this case the
singularities on the edges of triangles cannot be cancelled unless the
potential density is constant for all panels. 
Nonetheless one can use Galerkin BEM with constant panel basis functions for
hypersingular integral equations using condition (\ref{si7}). The rationale
for that is that the contributions related to the common edges should be
cancelled for linear or higher order methods, while the constant panel value
can be used as the magnitude for the parts which do not experience
cancellation.
\end{remark}

\subsubsection{Three-touch case}

As follows from the one- and two-touch cases the
three-touch case has the same singularity (\ref{si2}) as the two-touch case.
However, if we apply condition (\ref{si7}) the integrals are finite and due to (\ref{si3}) and (\ref{gdb8})%
\begin{equation}
M_{val}^{\prime }=2\sum_{j=1}^{3}l_{j}\ln \frac{p}{p-l_{j}}.  \label{si8}
\end{equation}%
Here we put subscript ``$val$'' to indicate that this is
not exactly the divergent integral, but the part computed under conditions (\ref{si7}). 

\section{Some numerical checks}

\subsection{Validation using Matlab library}

To verify the analytical expressions we implemented the entire algorithm in Matlab. The Matlab library function \emph{integral2} was used to evaluate 4D integrals (\ref{ps1})-(\ref{ps4}) directly. The tolerance for
the library function can be set by the user and for the cases reported below
are set to 0 for absolute tolerance and to 10$^{-12}$ for the relative
tolerance. We created several benchmark cases which realize different values
for the $\left\{ h_{1},...,h_{4}\right\} $ sets and provide the in Table \ref{table}.

\begin{table}[htbp]\label{table}
 \vspace*{-0.1in}
\begin{center}
\caption{Evaluation of no-touch Galerkin integrals}
 \vspace*{-0.1in}
\begin{tabular}{|l|l|l|l|}
	\hline 
Test case & 1 & 2 & 3 \\ 
\hline
$\mathbf{y}_{3}$ & $\left( 1/2,0,1+\sqrt{3}/2\right) $ & $\left( 1/2,\sqrt{6}%
/4,1+\sqrt{6}/4\right) $ & $\left( 1/2,-\sqrt{3}/2,1\right) $ \\ 
Cases called & 1,2,3,4 & 2,3,4,5 & 2,3,4,5,6,7 \\ 
$L$ & 0.139757030669707 & 0.149630247150535 & 0.156068357679434 \\ 
$M$ & 0.099860729206614 & 0.114715727210190 & 0.111863573921226 \\ 
$\mathbf{L}^{\prime }\cdot \mathbf{i}_{x}$ & 0 & 0 & 0 \\ 
$\mathbf{L}^{\prime }\cdot \mathbf{i}_{y}$ & 0.022035244796804 & 
0.010953212167802 & 0.055673013677787 \\ 
$\mathbf{L}^{\prime }\cdot \mathbf{i}_{z}$ & -0.099860729206614 & 
-0.114715727210190 & -0.111863573921226 \\ 
$M^{\prime }$ & 0.046564310284965 & 0.137859073743097 & -0.138417139905960
\\ 
$\epsilon _{\max }$ & $4.5\cdot 10^{-16}$ & $1.5\cdot 10^{-15}$ & $1.7\cdot
10^{-15}$\\
\hline
\end{tabular}
\end{center}
 \vspace*{-0.1in}
\end{table}

In these tests both triangles $S_{x}$ and $S_{y}$ are equilateral with the
side length $l_{xi}=l_{yj}=1$, $i,j=1,2,3.$ The vertices are $\mathbf{x}_{1}=
$ $\left( 0,0,0\right) $, $\mathbf{x}_{2}$ $=\left( 1,0,0\right) $, $\mathbf{%
x}_{3}$ $=$ $\left( 1/2,\sqrt{3}/2,0\right) $, $\mathbf{y}_{1}=\left(
1,0,1\right) $, $\mathbf{y}_{2}=$ $\left( 0,0,1\right) $, while $\mathbf{y}%
_{3}$ is specific for each case and given in Table 1. The value $\epsilon
_{\max }$ shows the absolute maximum difference between all integrals 
with the Matlab library. The ``Cases called''
show which cases from our classification were called to compute the values,
which we tried to select in the way to cover all possible realizations
\#1-\#7 for no-touch integrals. These errors are consistent
with the double precision accuracy.

\subsection{Touching cases}

One-, two, and three-touch cases were created to verify analytical
expressions and convergence for singular cases. The procedure was the
following. The theoretical or limiting value was computed using the present
method and it was compared with the corresponding no-touch case where the
difference between the touching vertices is $\epsilon $. In all cases the
source triangle is the equilateral triangle with vertices $\mathbf{x}_{1}=$ $%
\left( 0,0,0\right) $, $\mathbf{x}_{2}$ $=\left( 1,0,0\right) $, $\mathbf{x}%
_{3}$ $=$ $\left( 1/2,\sqrt{3}/2,0\right) $, while the receiver equilateral
triangle of the same size has vertices  $\mathbf{y}_{j}=$ $\mathbf{y}%
_{j0}+\epsilon \mathbf{i}_{z}$, $j=1,2,3,$ with $\mathbf{y}_{j0}$ specific
for the one-, two-, and three-touch cases. The relative difference is
computed as $\epsilon _{rel}=\left( F-F_{0}\right) /F_{0}$, $F=L,M,M^{\prime
},$ where $F_{0}$ is a limiting or asymptotic value.

\begin{figure}[tbh]
 \vspace*{-0.1in}
	\begin{center}
		\includegraphics[width=\textwidth,trim=0 2.6in 0 0]{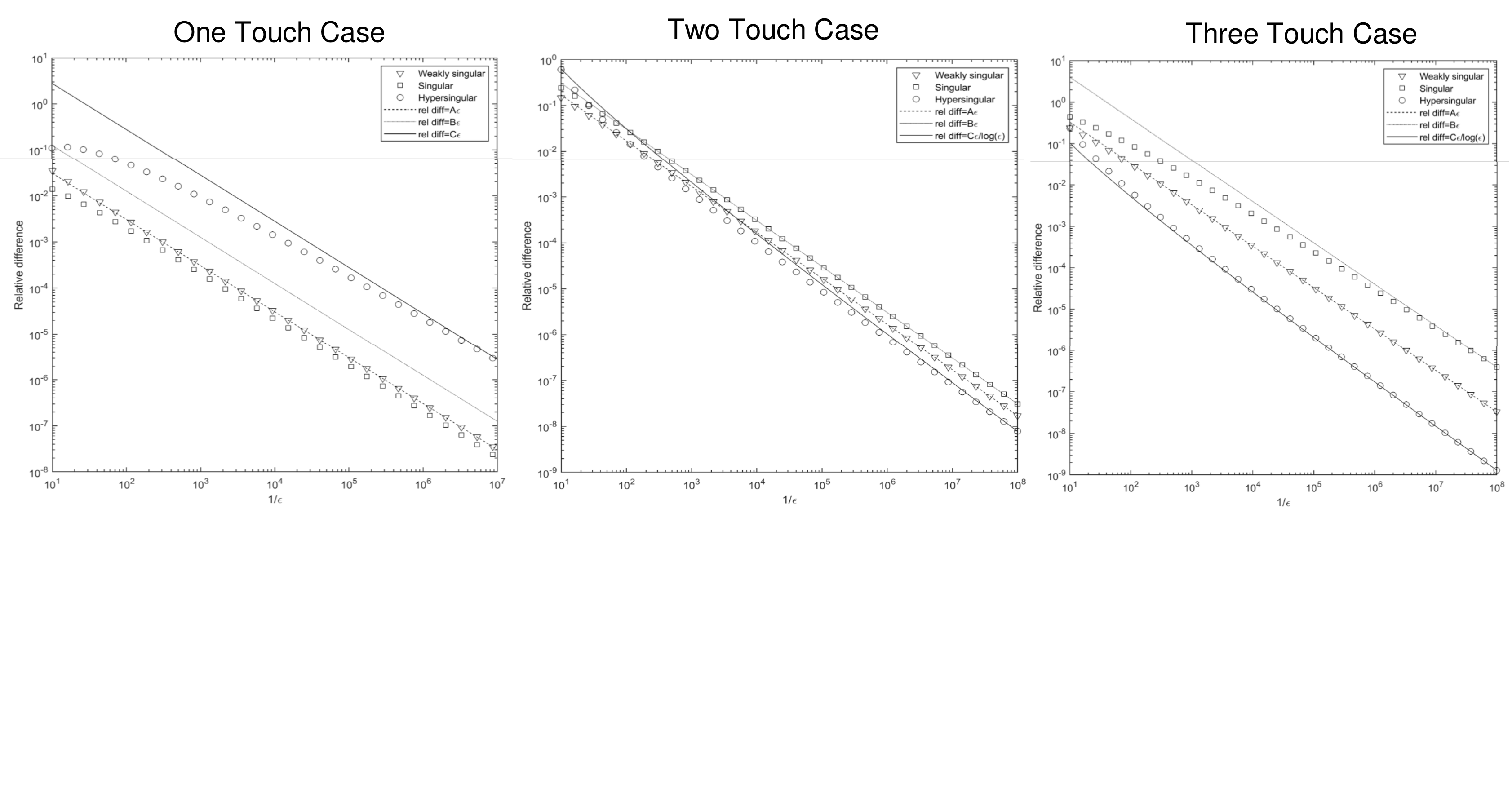}
	\end{center}
 \vspace*{-0.1in}
	\caption{Relative differences for nonzero $\epsilon$ and $\epsilon =0$ for $L$ (weakly singular), $M$ (singular), and $L^{\prime }$ (hypersingular) integrals for \\ Left - One Touch:  The lines show asymptotic dependences; \\ 
	Center - Two Touch: Asymptotics of the hypersingular integral is used  for the reference;\\ 
    Right - Three Touch: Self-action analytical and asymptotic values are used for all integrals as reference.}
	\label{Fig3}
\end{figure}

The convergence in the one-touch case for integrals $L,M$, and $M^{\prime }$
is illustrated in Fig. \ref{Fig3} (left). In this case $\mathbf{y}_{10}=\left(
0,0,0\right) ,$ $\mathbf{y}_{20}=\left( -1,0,0\right) ,$  $\mathbf{y}%
_{30}=\left( -1/2,0,\sqrt{3}/2\right) $. It is seen that all three integrals
converge to their limiting values $L_{0}=0.182526568122379$, $%
M_{0}=0.055671118815334$, $M_{0}^{\prime }=0.063116905873345$ computed at $%
\epsilon =0$ at linear convergence rate. The constants for the linear
convergence rates are computed using computed data for $\epsilon =10^{-7}$.

The convergence in the two-touch case for the same integrals is shown in
Fig. \ref{Fig3} (center). Here we set $\mathbf{y}_{10}=\left( 0,0,0\right) ,$ $%
\mathbf{y}_{20}=\left( 1,0,0\right) ,$ $\mathbf{y}_{30}=\left( 1/2,0,\sqrt{3}%
/2\right) $. The limiting values for the single and double layer potentials
computed at $\epsilon =0$ are $L_{0}=0.415922738854561$ and $%
M_{0}=0.706739910625218$. We also obtained $M_{0}^{\prime }=M_{00}^{\prime
}-H_{11}$, where $M_{00}^{\prime }=2.857471441252689$ is obtained at $%
\epsilon =0$ (see (\ref{gdb8}) under condition (\ref{si7}) for the common
edge), and $H_{11}=2\ln (2/\epsilon )+2$ comes from (\ref{si6}) for the
common edge length $l_{x1}=1$. Theoretically the convergence rates for $L$
and $M$ integrals are linear, while for $M^{\prime }$ this should be $%
\epsilon _{rel}\sim \epsilon /\ln (1/\epsilon )$. As in the previous case we
determined the constants at some value of $\epsilon $ ($\epsilon =10^{-8}$)
and plotted the curves on Fig. \ref{Fig3}. Again, the consistency of the computed
and theoretical rates is seen.

In the three-touch case illustrated in Fig. \ref{Fig3} (right) we have  $\mathbf{y}%
_{j0}=\mathbf{x}_{j}$, $j=1,2,3$. Here we have the following limiting
values, which can be computed straightforwardly (details are provided in the supplemental material), 
$L_{0}=\frac{3}{4}\ln 3$, and (\ref{si8}), (\ref{si6}), $M_{0}^{\prime
}=M_{val}^{\prime }-3H_{11}$, $M_{val}=6\ln 3$. Here factor 3 in front of $H_{11}$ is due to the contribution described by (\ref{si6}) as we have
three common edges in the three-touch case. It is remarkable that $M_{0}=0$,
which is a true self-action value for the double layer, cannot be used for
convergence study, as the double layer potential has a jump from $\pi $ on $S_{x}^{+}$ to $-\pi $ on $S_{x}^{-}$. Integration of this value over the
area of the receiver triangle provides $M_{0}=\pi A_{y}=\pi \sqrt{3}/2$ as a
limiting value as the receiver triangle approaches the source triangle from
the side denoted as $S_{x}^{+}$ (from positive $z$ for the oriented triangle
with the normal directed as $\mathbf{i}_{z}$). The expected convergence
rates in this case are the same as in the two-touch case. We determined the
constants by matching data at $\epsilon =10^{-8}$ and showed the
approximations along with the calculated data in Fig. \ref{Fig3}. 

\subsection{Comparison with \protect\cite{Lenoir2012SIAM}, \protect\cite%
{Salles2012matlab}}

Salles \cite{Salles2012matlab} provides a matlab code to compute integrals
calculated in \cite{Lenoir2012SIAM} for cases \#6 and \#7 (in our notation).
We executed this code for the cases from the ``example of use'' \cite%
{Salles2012matlab}. In these cases the vertices are $x_{1}=$ $\left(
0,0,0\right) $, $x_{2}$ $=\left( 0,1,0\right) $, $x_{3}$ $=$ $\left(
1,0,0\right) $, $y_{1}=\left( 0,0,h\right) $, $y_{2}=$ $\left( 0,1,h\right) $%
, $y_{3}=\left( -1,0,h\right) $ for $L$ integrals, while the receiver
triangle for $\mathbf{L}^{\prime }$ integral has vertices $y_{1}=\left(
-2,0.5,0.01\right) $, $y_{2}=$ $\left( -1,1,0.01\right) $, $y_{3}=\left(
-1,0,0.01\right) .$ In Table 2 we provide outputs of \cite{Salles2012matlab}
and ours in the same format.

\begin{table}[htbp]
\begin{center}
\begin{tabular}{|l|l|l|l|l|}
\hline
Type & h & \cite{Lenoir2012SIAM}, \cite{Salles2012matlab} & This work & $%
\epsilon $ \\ 
\hline
$L$ & 0 & 0.4154834934268203 &  0.4154834934268203 & 0 \\ 
$L$ & 10$^{-4}$ &  0.4154834087866362 & 0.4154834087866360 & $2.2\cdot
10^{-16}$ \\ 
$L$ & 10$^{-3}$ & 0.4154773308369880 &  0.4154773308369882 & $1.7\cdot
10^{-16}$ \\ 
$L$ & 10$^{-2}$ & 0.4150963397038615 & 0.4150963397038614 & $1.1\cdot
10^{-16}$ \\ 
$L$ & 10$^{-1}$ & 0.3986731498732934 & 0.3986731498732936 & $2.2\cdot
10^{-16}$ \\ 
$L$ & 1 &  0.1994877345160992 & 0.1994877345160997 & $5.6\cdot 10^{-16}$ \\ 
$\mathbf{L}^{\prime }\cdot \mathbf{i}_{x}$ & 10$^{-2}$ &  0.0937210251186380
& 0.0937210251186334 & $4.6\cdot 10^{-15}$ \\ 
$\mathbf{L}^{\prime }\cdot \mathbf{i}_{y}$ & 10$^{-2}$ & -0.0069668668015920
&  -0.0069668668016032 & $1.1\cdot 10^{-14}$ \\ 
$\mathbf{L}^{\prime }\cdot \mathbf{i}_{z}$ & 10$^{-2}$ & -0.0006289369951270
& -0.0006289369951278 & $7.8\cdot 10^{-16}$\\
\hline
\end{tabular}
\caption{Comparison of the results of the present work with some data in
literature \cite{Lenoir2012SIAM}, \cite{Salles2012matlab}. 
The absolute difference, $\epsilon$, between
the results of the present study and \cite{Salles2012matlab} for $L$ is of
the order of $10^{-16}$, i.e. consistent with double precision
arithmetics. The difference for $\mathbf{L}^{\prime }$ is 1-2 orders of
magnitude larger, but it is still small enough and can be explained by the
round-off errors in summation of several terms of different magnitude.}
\end{center}
\end{table}

\section{Conclusion}

This work shows that all integrals appearing in the Galerkin BEM
with constant basis functions for the Laplace Green's function can be evaluated analytically. 
However, one weakly singular integral formally can require 216 calls of the routine computing primitives. 
In practice, this can be reduced to 180 calls or lesser since some expansion coefficients are $s_{j0}$ are zero. 
Also, there exist some symmetries, which can be further investigated and used to improve efficiency. 
Furthermore, if simultaneous evaluation of single and double layer integrals and/or single layer derivatives is required a substantial 
amortization of computations occurs, as integrals in most cases can be computed simultaneously with little or no additional overhead.

For practical purposes the analytical integration can be employed for the nearfield integrals, 
while more efficient quadrature or expansion based \cite{Gumerov2021arXiv2} numerical integration can be applied for the far-field interactions. 
Nonetheless, the availability of the analytical methods to compute any kind of integrals is valuable for convergence studies and 
accurate evaluation of the order of the accuracy of the methods. In the future, the method can be 
extended to include linear and higher order basis function or include the Jacobians, 
which allow to account for curved panels, etc, as some work in this direction is already published (e.g., \cite{Salles2013Thesis}).

The method can be extended to other than Laplace Green's functions. Indeed, the recursive reduction process is the same independent of the kernel. 
All analytical challenges come only from construction of the PBF. A general expression for the PBF, however, is available (\ref{mp18}). 
In the case, when not all four integrals can be computed straightforwardly, first, some reduction may be still available, and, 
second, some special functions can be introduced to make the PBF computable analytically or semi-analytically. 
Recursive layer potential integrals via dimensionality reduction for the Laplace and Helmholtz kernel for higher order elements in the context of collocation BEM is discussed in our recent work~\cite{KGD2023RIPE}. 

\appendix

\section{Illustration of the algorithm applied to the computation of self-action integrals}

\subsection{Single layer potential}

Let us consider the case of self-action integral, when triangles $S_{x}$ and 
$S_{y}$ coincide. To verify the general solution we can also derive
analytical expressions for this case. Let us do this following
the algorithm.
\begin{itemize}
\item \textbf{Case definition}. Equation (\ref{s4}) shows that in the case%
\begin{equation}
\mathbf{X}_{0}=\mathbf{Y}_{0},\quad \mathbf{X}_{s}=\mathbf{Y}_{u},\quad 
\mathbf{X}_{t}=\mathbf{Y}_{v},  \label{sa1}
\end{equation}%
while (\ref{b01}) provides for the 4D case%
\begin{equation}
\mathbf{a}_{3} = \mathbf{-a}_{1}=-\mathbf{l}_{1},\quad \mathbf{a}_{4}=%
\mathbf{-a}_{2}=\mathbf{l}_{3},\quad \mathbf{e}_{4}=\mathbf{0},  \label{sa2}
\quad
U =I_{4}\left( \mathbf{l}_{1},-\mathbf{l}_{3},-\mathbf{l}_{1},\mathbf{l}%
_{3},\mathbf{0}\right) . 
\end{equation}%
Here we used (\ref{s11}) where the subscript $x$ (or $y$) is removed, since $%
S_{x}$ and $S_{y}$ coincide.

\item \textbf{4D to 3D reduction}. Since $\mathbf{e}_{4}=\mathbf{0}$ we have
from (\ref{mp3}) and (\ref{mp5}) 
\begin{equation}
s_{10}=s_{20}=s_{30}=s_{40}=h_{4}=0.  \label{sa3}
\end{equation}%
Relations (\ref{b2}) show then that only two 3D integrals contribute to the
overall integral $U$%
\begin{equation}\begin{aligned}
U=I_{32}+I_{35}&=I_{3}\left( \mathbf{l}_{1},\mathbf{-l}_{3},\mathbf{-l}_{2},\mathbf{-l}_{1}\right) +I_{3}\left( -\mathbf{l}_{1},\mathbf{l}_{3},\mathbf{l}_{2},\mathbf{l}_{1}\right) = 2I_{3}\left( \mathbf{l}_{1},\mathbf{-l}_{3},\mathbf{-l}_{2},\mathbf{-l}_{1}\right) .  \label{sa4}
\end{aligned}\end{equation}
Here, first, we used (\ref{b4}), (\ref{sa2}), and (\ref{s10}) to express the
arguments of $I_{3}$ integrals via $\mathbf{l}_{j}$ and, second, the fact
that according the definition (\ref{mp1}) and (\ref{mp2}) the integrands
depending on $R_{d}$ only (e.g. in (\ref{b3})) do not change if we change
the sign of all parameters to the opposite.

\item \textbf{3D to 2D reduction}. It is clear now, that $\mathbf{l}_{1}\in
span\left( \mathbf{l}_{1},\mathbf{l}_{2},\mathbf{l}_{3}\right) $, and so $%
h_{3}=0$. Furthermore, since $\mathbf{e}_{3}=-\mathbf{l}_{1}$, the
coefficients of expansion of $\mathbf{e}_{3}$ over vectors $\left( \mathbf{l}%
_{1},-\mathbf{l}_{2},-\mathbf{l}_{3}\right) $ can be selected as%
\begin{equation}
s_{10}=-1,\quad s_{20}=0,\quad s_{30}=0.  \label{sa5}
\end{equation}%
Relation (\ref{b5}) then shows that only two 2D integrals contribute to $%
I_{3}$, so%
\begin{equation}\begin{aligned}
U&=2I_{3}\left( \mathbf{l}_{1},\mathbf{-l}_{3},\mathbf{-l}_{2},\mathbf{-l}%
_{1}\right) = 2\left( I_{22}+I_{24}\right) = 2\left[ I_{2}^{(s)}\left( \mathbf{%
-l}_{3},\mathbf{-l}_{2},-\mathbf{l}_{1}\right) +I_{2}^{(t)}\left( \mathbf{l}%
_{1},\mathbf{-l}_{3},\mathbf{l}_{3}\right) \right] .  \label{sa6}
\end{aligned}\end{equation}%
where the arguments of the integrals are found using (\ref{b7}).

\item \textbf{2D to 1D reduction}. Again this clearly shows that since $%
\mathbf{l}_{1}=-\mathbf{l}_{2}-\mathbf{l}_{3}\in span\left( \mathbf{l}_{2},%
\mathbf{l}_{3}\right) $ and $\mathbf{l}_{2}\in span\left( \mathbf{l}_{1},%
\mathbf{l}_{2}\right) $ we have for both integrals in (\ref{sa6}) $h_{2}=0$.
Moreover, 
\begin{equation}
\text{for }I_{2}^{(s)}:s_{10}=-1,\quad s_{20}=-1,\quad \text{for }%
I_{2}^{(t)}:s_{10}=0,\quad s_{20}=-1.  \label{sa7}
\end{equation}%
From (\ref{sa6}), (\ref{b8}), and (\ref{b10}) we obtain%
\begin{equation}\begin{aligned}
U &= 2\left( I_{12}^{(s)}+I_{14}^{(s)}+I_{12}^{(t)}\right) =2\left[
I_{1}\left( \mathbf{-l}_{2},-\mathbf{l}_{1}\right) +I_{1}\left( \mathbf{-l}%
_{3},-\mathbf{l}_{1}\right) +I_{1}\left( \mathbf{l}_{1},\mathbf{l}%
_{3}\right) \right]  \label{sa8} \\
&= 2\sum_{j=1}^{3}I_{1}\left( \mathbf{l}_{j},\mathbf{e}_{1}^{(j)}\right)
,\quad \mathbf{e}_{1}^{(1)}=\mathbf{l}_{3},\quad \mathbf{e}_{1}^{(2)}=%
\mathbf{e}_{1}^{(3)}=\mathbf{l}_{1}.  
\end{aligned}\end{equation}%
Here for convenience we rearranged summation and changed the signs of the
arguments in two integrals (see the comment to (\ref{sa4})).

\item \textbf{Evaluation}. First, we should determine the parallel and
perpendicular projections of vectors $\mathbf{e}_{1}^{(j)}$. Generally, for
projections we have (see Appendix)%
\begin{equation}\begin{aligned}
&\mathbf{e}_{1\Vert }^{(j)}=s_{10}^{(j)}\mathbf{l}_{j},\quad
h_{1}^{(j)}=\left| \mathbf{e}_{1\bot }^{(j)}\right| =\left| \mathbf{e}%
_{1}^{(j)}-\mathbf{e}_{1\Vert }^{(j)}\right| ,\quad s_{10}^{(j)}=\frac{%
\mathbf{e}_{1}^{(j)}\cdot \mathbf{l}_{j}}{l_{j}^{2}},\quad l_{j}=\left| 
\mathbf{l}_{j}\right| ,\quad j=1,2,3.  \label{sa9}
\end{aligned}\end{equation}%
Since for all integrals $h_{2}=h_{3}=h_{4}=0$, we have from the first case
of (\ref{k6}) and (\ref{b11}) 
\begin{equation}\begin{aligned}
&I_{1}\left( \mathbf{l}_{j},\mathbf{e}_{1}^{(j)}\right) = \left(
1+s_{10}^{(j)}\right) F_{1}\left( P_{1}^{(j)}\left( 1\right) \right)
-s_{10}^{(j)}F_{1}\left( P_{1}^{(j)}\left( 0\right) \right) ,\\
&P_{1}^{(j)}\left( s_{1}\right) =l_{j}\left|s_{1}+s_{10}^{(j)}\right|,
\label{sa10} \quad F_{1}\left( P_{1}^{(j)}\left( s_{1}\right) \right) = \frac{1}{12P_{1}^{(j)}%
}\ln \frac{\left( P_{1}^{(j)}+R_{1}^{(j)}\right) ^{2}}{h_{1}^{\left(
j\right) 2}},\\
&R_{1}^{(j)}\left( s_{1}\right) =\sqrt{P_{1}^{(j)2}\left( s_{1}\right) +h_{1}^{\left( j\right) 2}}, 
\quad j=1,2,3.
\end{aligned}\end{equation}

\item \textbf{Simplification}. Some simplifications can be
obtained for the integrals listed above. We have
\begin{equation}
F_{1}\left( P_{1}^{(j)}\right) =\frac{1}{12P_{1}^{(j)}}\ln \frac{\left(
R_{1}^{(j)}+P_{1}^{(j)}\right) ^{2}}{R_{1}^{(j)2}-P_{1}^{(j)2}}=\frac{1}{%
12P_{1}^{(j)}}\ln \frac{R_{1}^{(j)}+P_{1}^{(j)}}{R_{1}^{(j)}-P_{1}^{(j)}}.
\label{sa11}
\end{equation}%
This expression is invariant with respect to the sign of $P_{1}^{(j)}$, and
so instead of $P_{1}^{(j)}\left( s_{1}\right) $ we can use the following
signed function decorated by tilde%
\begin{equation}
\widetilde{P}_{1}^{(j)}\left( s_{1}\right) =l_{j}\left(
s_{1}+s_{10}^{(j)}\right) ,\quad F_{1}\left( P_{1}^{(j)}\right) =F_{1}\left( 
\widetilde{P}_{1}^{(j)}\right) =\frac{1}{12\widetilde{P}_{1}^{(j)}}\ln \frac{%
R_{1}^{(j)}+\widetilde{P}_{1}^{(j)}}{R_{1}^{(j)}-\widetilde{P}_{1}^{(j)}}.
\label{sa12}
\end{equation}%
Hence, integral in (\ref{sa10}) turns into 
\begin{equation}\begin{aligned}
I_{1}\left( \mathbf{l}_{j},\mathbf{e}_{1}^{(j)}\right) =&\left(
1+s_{10}^{(j)}\right) F_{1}\left( \widetilde{P}_{1}^{(j)}\left( 1\right)
\right) -s_{10}^{(j)}F_{1}\left( \widetilde{P}_{1}^{(j)}\left( 0\right)
\right)  \label{sa14} \\
=&\frac{1}{12l_{j}}\ln \frac{R_{1}^{(j)}\left( 1\right) +\widetilde{P}%
_{1}^{(j)}\left( 1\right) }{R_{1}^{(j)}\left( 1\right) -\widetilde{P}%
_{1}^{(j)}\left( 1\right) }\frac{R_{1}^{(j)}\left( 0\right) -\widetilde{P}%
_{1}^{(j)}\left( 0\right) }{R_{1}^{(j)}\left( 0\right) +\widetilde{P}%
_{1}^{(j)}\left( 0\right) }. 
\end{aligned}\end{equation}%
It follows from (\ref{sa8}) and (\ref{sa9}) that%
\begin{equation}\begin{aligned}
&s_{10}^{(1)}=\frac{\mathbf{l}_{1}\cdot \mathbf{l}_{3}}{l_{1}^{2}}=\frac{l_{2}^{2}-l_{1}^{2}-l_{3}^{2}}{2l_{1}^{2}},\quad s_{10}^{(2)}=\frac{\mathbf{l}_{1}\cdot \mathbf{l}_{2}}{l_{2}^{2}}=\frac{l_{3}^{2}-l_{1}^{2}-l_{2}^{2}}{2l_{2}^{2}},\quad
s_{10}^{(3)}=\frac{\mathbf{l}_{1}\cdot \mathbf{l}_{3}}{l_{3}^{2}}=\frac{l_{2}^{2}-l_{1}^{2}-l_{3}^{2}}{2l_{3}^{2}},  \label{sa15}
\end{aligned}\end{equation}%
where we got rid of scalar products using the cosine theorem for the
triangle (or $\left| \mathbf{l}_{1}+\mathbf{l}_{2}\right| ^{2}=\left| 
\mathbf{l}_{3}\right| ^{2}$) 
\begin{equation}
l_{1}^{2}+2\mathbf{l}_{1}\cdot \mathbf{l}_{2}+l_{2}^{2}=l_{3}^{2},\quad 
\mathbf{l}_{1}\cdot \mathbf{l}_{2}=\frac{l_{3}^{2}-l_{1}^{2}-l_{2}^{2}}{2}.
\label{sa16}
\end{equation}%
Furthermore, we have the following expressions for functions $\widetilde{P}%
_{1}^{(j)}$ and $R^{(j)}$ 
\begin{equation}\begin{aligned}
R_{1}^{(1)}\left( 0\right) =& l_{3},\quad R_{1}^{(2)}\left( 0\right)
=l_{1},\quad R_{1}^{(3)}\left( 0\right) =l_{1},  \quad
R_{1}^{(1)}\left( 1\right) = l_{2},\quad R_{1}^{(2)}\left( 1\right)
=l_{3},\\ R_{1}^{(3)}\left( 1\right) &=l_{2},\quad
\widetilde{P}_{1}^{(j)}\left( 0\right) = l_{j}s_{j0}^{(1)},\quad \widetilde{%
P}_{1}^{(j)}\left( 1\right) =l_{j}\left( 1+s_{10}^{(j)}\right) ,\quad
j=1,2,3.\quad  \label{sa17}
\end{aligned}\end{equation}%
Substituting these into expressions for integrals (\ref{sa14}) and
simplifying, we obtain
\begin{equation}
I_{1}\left( \mathbf{l}_{j},\mathbf{e}_{1}^{(j)}\right) =\frac{1}{6l_{j}}\ln 
\frac{l_{1}+l_{2}+l_{3}}{l_{1}+l_{2}+l_{3}-2l_{j}},\quad j=1,2,3,
\label{sa18}
\end{equation}%
The final expressions (\ref{sa8}) and (\ref{s6}), therefore, can be written
in the form%
\begin{equation}
L=4A^{2}U,\quad U=\frac{1}{3}\sum_{j=1}^{3}\frac{1}{l_{j}}\ln \frac{p}{%
p-l_{j}},\quad p=\frac{1}{2}\sum_{j=1}^{3}l_{j},  \label{sa19}
\end{equation}%
where $p$ is the half of the perimeter and $A$ is the area of triangle $%
S_{x}.$ This is consistent with the expression derived in \cite%
{Eibert1995IEEE} and simplified in \cite{Sievers2005IEEE} (the expression
under logs there can be slightly simplified to obtain (\ref{sa19}; also the
modulus under log can be removed, since the expression is always positive
due to the triangle inequality $p>l_{j},$ $j=1,2,3$).
\end{itemize}

\subsection{Gradient of single layer and double layer potential}

Similarly to the single layer case we can consider for illustration
computation of the self-action integral. As it is mentioned in this case we
have $M=0$. This case is characterized by $h_{4}=h_{3}=h_{2}=0$. So, based
on $\mathbf{a}_{1},...,\mathbf{a}_{4},\mathbf{e}_{4}$ from (\ref{sa2}) we
can determine parameters of 3D integrals, which are%
\begin{equation}
F_{1} = 6l_{1}AI_{3}\!\left( \mathbf{-l}_{1},\mathbf{l}_{3},\mathbf{l}_{1},%
\mathbf{0}\right) ,  \label{db10} \quad\!\!\!
F_{2} = 6l_{2}AI_{3}\!\left( \mathbf{-l}_{1},\mathbf{l}_{3},\mathbf{l}_{2},%
\mathbf{l}_{1}\right) ,  \quad\!\!\!
F_{3} =6l_{3}AI_{3}\!\left( \mathbf{-l}_{1},\mathbf{l}_{3},\mathbf{-l}_{3},%
\mathbf{0}\right) , 
\end{equation}%
where we used relation (\ref{db7.3}) expressing $\widehat{I}_{3}$ via $I_{3}$
and removed subscript $x$ as $S_{x}$ and $S_{y}$ coincide. Following the
same procedure as (\ref{sa5})-(\ref{sa19}) we can find that%
\begin{equation}
F_{i}=\frac{3l_{i}}{4A}L,\quad i=1,2,3,  \label{db11}
\end{equation}%
where $L$ is from (\ref{sa19}). Final summation then yields%
\begin{equation}
\mathbf{F}=\sum_{i=1}^{3}\mathbf{n}_{ci}F_{i}=-\mathbf{n\times }%
\sum_{i=1}^{3}\frac{\mathbf{l}_{i}}{l_{i}}F_{i}=\frac{3L}{4A}\mathbf{n}%
\times \sum_{i=1}^{3}\mathbf{l}_{i}=\mathbf{0.}  \label{db12}
\end{equation}

\section{Method for identifying the parameters}

Given some vector $\mathbf{e}_{d}\in \mathbb{R}^{3}$ and a set of vectors $%
\left\{ \mathbf{a}_{1},...,\mathbf{a}_{d}\right\} $, $\mathbf{a}_{j}\in 
\mathbb{R}^{3}$, $\mathbf{a}_{j}\neq 0,$ $j=1,...,d,$ we need to determine
the perpendicular and parallel components $\mathbf{e}_{d\bot }$ and $\mathbf{%
e}_{d||}$, $\mathbf{e}_{d}=\mathbf{e}_{d||}+\mathbf{e}_{d\perp }$, such that 
$\mathbf{e}_{d||}\in span\left( \mathbf{a}_{1},...,\mathbf{a}_{d}\right) $,
while $\mathbf{e}_{d\perp }\cdot \mathbf{a}_{j}=0$. Moreover, we need to
find the coefficients of expansion $s_{10},...,s_{d0}$ for the parallel
component $\mathbf{e}_{d||}=s_{10}\mathbf{a}%
_{1}+...+s_{d0}\mathbf{a}_{d}$. 
$h_{d}$ is defined as the norm of
vector $\mathbf{e}_{d\perp }$. We briefly describe how to do this using Gram-Schmidt
orthogonalization.

The original Gram-Schmidt algorithm is designed
for orthogonalization of a linearly independent vector set, while here
 some vectors $\mathbf{a}_{1},...,\mathbf{a}_{d}$ are linearly 
dependent.
So, the proposed procedure consist of two steps: 
orthogonalization, and, determination of the expansion coefficients.

\subsection{Orthogonalization}

The standard Gram-Schmidt process is given as%
\begin{equation}
\mathbf{u}_{1}=\mathbf{a}_{1},\quad \mathbf{u}_{i}=\mathbf{a}%
_{i}-\sum_{k=1}^{i-1}c_{ki}\mathbf{u}_{k},\quad i=2,...,d,\quad c_{ki}=\frac{%
\mathbf{u}_{k}\cdot \mathbf{a}_{i}}{\left| \mathbf{u}_{k}\right| ^{2}},
\label{ap1}
\end{equation}%
where $\left\{ \mathbf{u}_{1},...,\mathbf{u}_{d}\right\} $ is the orthogonal
set, such as $span\left( \mathbf{u}_{1},...,\mathbf{u}_{d}\right)
=span\left( \mathbf{a}_{1},...,\mathbf{a}_{d}\right) $ and $c_{ki}\mathbf{u}%
_{k}$ is the projection of $\mathbf{a}$ on $\mathbf{u}_{k}$. As the steps $%
i=2,...,d$ are applied recursively all $\left\{ \mathbf{u}_{1},...,\mathbf{u}%
_{d}\right\} $ can be obtained for the linearly independent set $\left\{ 
\mathbf{a}_{1},...,\mathbf{a}_{d}\right\} $. In case if $\mathbf{a}_{j}$
depends linearly on $\mathbf{a}_{1},...,\mathbf{a}_{j-1}$(i.e., $\mathbf{a}%
_{j}\in span\left\{ \mathbf{a}_{1},...,\mathbf{a}_{j-1}\right\} $) the
process produces $\mathbf{u}_{j}=\mathbf{0}$ and cannot be continued as
represented by (\ref{ap1}), since $\left| \mathbf{u}_{j}\right| =0$. 
An easy way to modify the process which can pass through linearly dependent vectors
is to redefine the projection coefficients as%
\begin{equation}
c_{ki}=\frac{\mathbf{u}_{k}\cdot \mathbf{a}_{i}}{u_{k}^{2}},\quad
u_{k}^{2}=\left\{ 
\begin{array}{c}
\left| \mathbf{u}_{k}\right| ^{2},\quad \left| \mathbf{u}_{k}\right| \neq 0
\\ 
1,\quad \left| \mathbf{u}_{k}\right| =0%
\end{array}%
\right. ,\quad i,k=1,...,d  \label{ap2}
\end{equation}%
The process (\ref{ap1}) produces the orthogonal set $%
\left\{ \mathbf{u}_{1},...,\mathbf{u}_{d}\right\} $ with $\mathbf{u}_{j}$ $=%
\mathbf{0}$ in place of $\mathbf{a}_{j}\in span\left\{ \mathbf{a}_{1},...,%
\mathbf{a}_{j-1}\right\} .$

\subsection{Expansion coefficients}

We can now expand all vectors $\mathbf{a}_{1},...,\mathbf{a}%
_{d}$ over the set $\left\{ \mathbf{u}_{1},...,\mathbf{u}_{d}\right\} $ as 
\begin{equation}
\mathbf{a}_{i}=\sum_{k=1}^{i}c_{ki}\mathbf{u}_{k},\quad i=1,...,d, 
\label{ap3}
\end{equation}%
considering the expansion%
\begin{equation}
\mathbf{e}_{d}=\sum_{i=1}^{d}s_{i0}\mathbf{a}_{i}+\mathbf{e}_{d\bot
}=\sum_{i=1}^{d}s_{i0}\sum_{k=1}^{i}c_{ki}\mathbf{u}_{k}+\mathbf{e}_{d\bot }.
\label{ap4}
\end{equation}%
Taking the scalar product of $\mathbf{e}_{d}$ and $\mathbf{u}_{j}$ we obtain%
\begin{equation}
\sum_{i=1}^{d}A_{ji}s_{i0}=b_{j},\quad A_{ji}=c_{ji}\left| \mathbf{u}%
_{j}\right| ^{2},\quad b_{j}=\mathbf{u}_{j}\cdot \mathbf{e}_{d},\quad
j=1,...,d.  \label{ap5}
\end{equation}%
This is a linear system for unknowns $s_{i0}$ with the upper triangular
matrix with entries $A_{ji}$ and the right hand side vector with components $%
b_{j}$. It can be solved recursively. For the $j$th row 
corresponding to $\mathbf{u}_{j}=\mathbf{0}$ we have $A_{ji}=0$ and $b_{j}=0$%
; \ so, one should set $s_{j0}=0$. The solution is then 
\begin{equation}
s_{j0}=\left\{ 
\begin{array}{c}
\frac{1}{A_{jj}}\left( b_{j}-\sum_{i=j+1}^{d}A_{ji}s_{i0}\right) ,\quad
A_{jj}\neq 0 \\ 
0,\quad A_{jj}=0%
\end{array}%
\right. ,\quad j=d,d-1,...,1.  \label{ap6}
\end{equation}%
As soon as all $s_{i0}$ are determined, parameter $h_{d}$ can be found from%
\begin{equation}
h_{d}=\left| \mathbf{e}_{d\bot }\right| =\left| \mathbf{e}%
_{d}-\sum_{i=1}^{d}s_{i0}\mathbf{a}_{i}\right| .  \label{ap7}
\end{equation}

\end{document}